\documentclass[12pt,a4paper]{amsart}
\usepackage[utf8]{inputenc}
\usepackage[foot]{amsaddr}

\usepackage{mathrsfs,amsthm,xcolor,verbatim,bbm,amsmath,amsfonts,
amssymb,nicefrac,enumitem,hyperref,bm,mathtools,xparse,etoolbox}
\usepackage[margin=1in]{geometry}
\usepackage[capitalise,sort]{cleveref} 

\crefname{enumi}{item}{items}

\crefname{equation}{}{}
\crefname{subsection}{Subsection}{Subsections}
\crefname{prop}{Proposition}{Propositions}

\hypersetup{
    colorlinks,
    linkcolor={red!50!black},
    citecolor={green},
    urlcolor={blue!80!black}
}

\theoremstyle{plain}
\newtheorem{theorem}{Theorem} [section]
\newtheorem{lemma}[theorem]{Lemma}
\newtheorem{prop}[theorem]{Proposition}
\newtheorem{cor}[theorem]{Corollary}

\theoremstyle{definition}
\newtheorem{definition}[theorem]{Definition}
\newtheorem{case}{Case}
\newtheorem{remark}[theorem]{Remark}
\newtheorem{example}{Example}[section]

\crefname{case}{Case}{Cases}

\numberwithin{equation}{section}

\DeclareFontEncoding{LS1}{}{}
\DeclareFontSubstitution{LS1}{stix}{m}{n}
\DeclareMathAlphabet{\mathscr}{LS1}{stixscr}{m}{n}

\newcommand{\R}{\mathbb{R}}
\newcommand{\N}{\mathbb{N}}

\newcommand{\ssuml}{\textstyle\sum\limits}
\newcommand{\sprod}{\textstyle\prod}
\newcommand{\with}{\curvearrowleft}

\newcommand{\cA}{\mathcal{A}}
\newcommand{\cI}{\mathcal{I}}
\newcommand{\cO}{\mathcal{O}}
\newcommand{\cP}{\mathcal{P}}
\newcommand{\cR}{\mathcal{R}}

\newcommand{\bfk}{\mathbf{k}}

\newcommand{\bfI}{\mathbf{I}}
\newcommand{\bfM}{\mathbf{M}}
\newcommand{\bfN}{\mathbf{N}}
\newcommand{\bfP}{\mathbf{P}}

\newcommand{\scrP}{\mathscr{P}}

\newcommand{\fC}{\mathfrak{C}}
\newcommand{\fL}{\mathfrak{L}}

\newcommand{\fr}{\mathfrak{r}}
\newcommand{\fd}{\mathfrak{d}}
\newcommand{\fm}{\mathfrak{m}}
\newcommand{\fp}{\mathfrak{p}}

\DeclarePairedDelimiter{\norm}{\lVert}{\rVert}
\DeclarePairedDelimiter{\abs}{\lvert}{\rvert}
\DeclarePairedDelimiter{\rbr}{(}{)}
\DeclarePairedDelimiter{\br}{[}{]}
\DeclarePairedDelimiter{\cu}{\{}{\}}

\newcommand{\qqandqq}{\qquad\text{and}\qquad}

\newcommand{\num}[1]{\br*{ #1 } }

\newcommand{\cost}[1]{\operatorname{Cost}_{#1}}

\ExplSyntaxOn

\bool_new:N \g_noteobserve

\NewDocumentCommand{\setnote}{}{
  \bool_gset_true:N \g_noteobserve
}

\NewDocumentCommand{\setobserve}{}{
  \bool_gset_false:N \g_noteobserve
}

\NewDocumentCommand{\nobs}{ o }{
  \IfValueT{#1}{
    \str_if_eq:noTF {note} {#1} {
      \bool_gset_true:N \g_noteobserve
    } {
      \str_if_eq:noTF {Note} {#1} {
        \bool_gset_true:N \g_noteobserve
      } {
        \bool_gset_false:N \g_noteobserve
      }
    }
  }
  \bool_if:nTF { \g_noteobserve } {
    \bool_gset_false:N \g_noteobserve
    note
  } {
    \bool_gset_true:N \g_noteobserve
    observe
  }
  \IfValueF{#1}{~}
}

\NewDocumentCommand{\Nobs}{ o }{
  \IfValueT{#1}{
    \str_if_eq:noTF {note} {#1} {
      \bool_gset_true:N \g_noteobserve
    } {
      \str_if_eq:noTF {Note} {#1} {
        \bool_gset_true:N \g_noteobserve
      } {
        \bool_gset_false:N \g_noteobserve
      }
    }
  }
  \bool_if:nTF { \g_noteobserve } {
    \bool_gset_false:N \g_noteobserve
    Note
  } {
    \bool_gset_true:N \g_noteobserve
    Observe
  }
  \IfValueF{#1}{~}
}

\int_new:N \g_furthermore

\NewDocumentCommand{\Moreover}{ o o }{
  \IfValueT{#1}{
    \str_case:nn {#1} {
      {Furthermore} {\int_set:Nn {\g_furthermore} {0}}
      {Moreover} {\int_set:Nn {\g_furthermore} {1}}
      {In~addition} {\int_set:Nn {\g_furthermore} {2}}
      {note} {\bool_gset_true:N \g_noteobserve}
      {observe} {\bool_gset_false:N \g_noteobserve}
    }
    \IfValueT{#2}{
      \str_case:nn {#2} {
        {Furthermore} {\int_set:Nn {\g_furthermore} {0}}
        {Moreover} {\int_set:Nn {\g_furthermore} {1}}
        {In~addition} {\int_set:Nn {\g_furthermore} {2}}
        {note} {\bool_gset_true:N \g_noteobserve}
        {observe} {\bool_gset_false:N \g_noteobserve}
      }
    }
  }
  \int_case:nn { \int_mod:nn {\g_furthermore} {3} } {
    { 0 } { Furthermore,~\nobs that}
    { 1 } { Moreover,~\nobs that}
    { 2 } { In~addition,~\nobs that}
  }
  \int_incr:N \g_furthermore
  \IfValueF{#1}{~}
}

\bool_new:N \g_hencetherefore

\NewDocumentCommand{\hence}{}{
  \bool_if:nTF { \g_hencetherefore } {
    \bool_gset_false:N \g_hencetherefore
    hence~
  } {
    \bool_gset_true:N \g_hencetherefore
    therefore~
  }
}

\NewDocumentCommand{\Hence}{}{
  \bool_if:nTF { \g_hencetherefore } {
    \bool_gset_false:N \g_hencetherefore
    Hence,~we~obtain~
  } {
    \bool_gset_true:N \g_hencetherefore
    Therefore,~we~obtain~
  }
}

\ExplSyntaxOff

\NewDocumentEnvironment{cproof}{m}
{\begin{proof}[Proof of \cref{#1}]}%
	{\noindent The proof of \cref{#1} is thus complete.
\end{proof}}

\NewDocumentEnvironment{cproof2}{m}
{\begin{proof}[Proof of \cref{#1}]}%
	{\noindent This completes the proof of \cref{#1}.
\end{proof}}

\title[DNN approximation of composite functions]{Deep neural network approximation of composite functions without the curse of dimensionality}

\author{Adrian Riekert}

\email{ariekert@uni-muenster.de}

\address{Applied Mathematics: Institute for Analysis and Numerics, University of M{\"u}nster, Germany}

\date{\today}

\begin{document}
	
	\begin{abstract}
		In this article we identify a general class of high-dimensional continuous functions that can be approximated by deep neural networks (DNNs) with the rectified linear unit (ReLU) activation
		without the curse of dimensionality.
		In other words, the number of DNN parameters grows at most polynomially in the input dimension and the approximation error. 
		The functions in our class can be expressed as a potentially unbounded number of compositions of special functions which include products, maxima, and certain parallelized Lipschitz continuous functions.
	\end{abstract}

\maketitle


\section{Introduction}
Many practically relevant numerical approximation algorithms suffer from the curse of dimensionality (cf., e.g., Bellman~\cite{Bellman1957} and Novak \& Ritter~\cite{NovakRitter1997}).
Roughly speaking, this means that the number of parameters needed to approximate certain functions grows exponentially in the input dimension, which is often problematic in high dimensions.
In recent years deep neural networks (DNNs)
have been successfully employed in numerous high-dimensional applications, where the unknown function to be approximated can depend on thousands or even millions of real parameters. It thus seems that DNNs are able to overcome the curse of dimensionality in many relevant cases.
Nevertheless, the reasons for these promising practical results are still not fully understood.

In this article we aim to enhance the understanding of the approximation capabilities of DNNs by identifying a suitably large class of continuous functions that can provably be approximated by DNNs without the curse of dimensionality. That is, the number of DNN parameters is allowed to grow at most polynomially with respect to the input dimension and the prescribed approximation accuracy. 
The functions in this class can be expressed as compositions of particular functions which include products, maxima, and parallelized Lipschitz continuous functions; and the number of functions in the composition itself is allowed to grow polynomially in the input dimension, which leads to interesting new examples of approximable functions.

\subsection{Literature overview}
The fact that neural networks with a suitable activation function and sufficiently many parameters can approximate any continuous function up to an arbitrarily small error is known in the literature as the universal approximation theorem. Qualitative results of this type were first established, e.g., in \cite{Cybenko1989, FUNAHASHI1989183, HORNIK1991251, HORNIK1989359, LESHNO1993861}.
We also refer to Pinkus~\cite{Pinkus1999Approx} for further references regarding such universal approximation results.

Quantitative upper bounds on the number of DNN parameters needed to approximate continuous functions can be found, e.g., in  \cite{GuehringRaslanKutyniok2020,ZhouPuWang2017NIPS,Pinkus1999Approx,ShenYangZhang2020,Yarotsky2018optimal}.
These estimates suffer from the curse of dimensionality, since one needs $\Omega ( \varepsilon ^{ - c d } ) $ DNN parameters to approximate a general Lipschitz function $f \in C ([ 0 , 1 ] ^d , \R )$ with accuracy $\varepsilon> 0$, for some constant $c > 0$. Without further assumptions on $f$, this is essentially unavoidable since these approximation rates are optimal for all commonly used activation functions, as shown, e.g., in \cite{Mhaskar1996, Pinkus1999Approx} for the case of shallow networks (i.e., networks with one hidden layer)
 and, e.g., in \cite{ShenYangZhang2020,YAROTSKY2017} for the case of deep networks.
Thus, in order to overcome the curse of dimensionality one needs to restrict to special continuous functions with additional properties.

One possibility is to assume certain regularity properties. Indeed, for a sufficiently smooth function $f \in C^s ( [ 0 , 1 ] ^d , \R )$ with Sobolev norm bounded by a constant it can be shown (cf., e.g., \cite{Guehring2020Approx,JiaoShenLinHuang2021,LuShenYangZhang2021,Mhaskar1996,MontanelliDu2019,PETERSEN2018,YAROTSKY2017}) that $\cO ( \varepsilon ^{ - d / 2 s } )$ parameters are sufficient to achieve an accuracy of $\varepsilon$, up to a prefactor independent of $\varepsilon$ and the target function $f$ (that does, however, grow exponentially in $d$). Thus if $f$ is sufficiently regular, e.g.~if $s \ge d/2 $, one can obtain a polynomial approximation rate in terms of $\varepsilon $.

Another idea is to focus on functions with a certain hierarchical structure, which can be exploited by suitable DNN architectures.
For approximation results in the context of functions that can be written as an iterated composition of certain two-dimensional functions and thus be represented by a binary tree structure we refer, e.g., to \cite{MhaskarLiaoPoggio2016,Mhaskar_Liao_Poggio_2017,MhaskarPoggio2016,Poggio2017WhyAW}.
Similarly, approximation results for suitable compositional functions represented by directed graphs have been established in \cite{Gong2023approximation,KangGong2022}.
 The ideas we use in this article are in some sense related to this approach.

The fact that in particular deep network architectures have advantageous approximation properties compared to shallow networks has also been established in different settings, e.g., in \cite{Daniely2017,Eldan2016,GraeberJentzen2023,IbragimovKoppensteiner2023,SafranShamir2017}.

Finally, we refer to Gühring et al.~\cite{GuehringRaslanKutyniok2020} for a more extensive overview and further references on the approximation properties of DNNs.

\subsection{Notation}
\label{subsec:dnn:notation}
In this article we consider fully connected feedforward deep neural networks (DNNs) with the multidimensional \emph{rectified linear unit} (ReLU) activation function $\fr \colon \bigcup_{d \in \N} \R^d \to \bigcup_{d \in \N} \R^d$
which satisfies for all 
$d \in \N$, $x = (x_1, \ldots, x_d ) \in \R^d$
that
$\fr ( x ) = (\max \cu{x_1, 0 }, \ldots, \allowbreak \max \cu{x_d , 0 }) $.
The following mathematical description of DNNs,
which we employ throughout the remainder of this article, is inspired by Petersen \& Voigtlaender~\cite[Definition 2.1]{PETERSEN2018}
(see also, e.g., Beck et al.~\cite[Definition 2.9]{BeckJentzenKuckuck2022}).
We denote by
$\bfN$ the set of DNNs, which is given by
\begin{equation}
\textstyle
    \bfN = \bigcup_{L \in \N} \bigcup_{\ell_0, \ell_1, \ldots, \ell_L \in \N} \rbr*{ \bigtimes_{k=1}^{L} \rbr*{ \R^{\ell_k \times \ell_{k-1}} \times \R^{\ell_k} } } .
\end{equation}
For a neural network
\begin{equation*}
    \Phi = ((W_1,b_1), (W_2,b_2), \ldots,  (W_L, b_L)) \in \textstyle \rbr*{ \bigtimes_{k=1}^{L} \rbr*{ \R^{\ell_k \times \ell_{k-1}} \times \R^{\ell_k} } }
\end{equation*}
we define
its input dimension $\cI(\Phi) = \ell_0 \in \N$,
its output dimension
$\cO(\Phi) = \ell_L \in \N$,
and
 its layer functions
$\cA_i^\Phi \colon \R^{ \ell_{i-1} } \to \R^{ \ell_i }$,
$i \in \{ 1, 2, \ldots, L \}$,
by $\cA_i^\Phi ( x ) = W_i x + b_i$.
The \emph{realization} of $\Phi$ is the function $\cR ( \Phi ) \in C ( \R^{ \cI ( \Phi ) } , \R^{ \cO ( \Phi ) } )$
defined as the composition of the affine linear layer functions and the nonlinear activation $\fr$,
and thus
satisfies for all
$x \in \R^{ \cI ( \Phi ) }$
that
\begin{equation*}
\cR ( \Phi ) ( x ) = \cA^\Phi_L ( \fr ( \cA^\Phi_{L-1} ( \fr ( \cdots ( \fr ( \cA^\Phi_1 ( x ) ) ) \cdots ) ) ) ) \in \R ^{ \cO ( \Phi ) } .
\end{equation*}
The real entries of the matrices $W_i$ are also called \emph{weights} and the entries of the vectors $b_i$ are called \emph{biases}. 
The number of parameters of $\Phi$,
which equals the overall number of real numbers describing the weights and biases,
is
\begin{equation*}
\cP(\Phi) = \ssuml_{i=1}^L \ell_i( \ell_{i-1}+1) \in \N .
\end{equation*}
This is the measure we use to describe the complexity of the neural network $\Phi$.

For every $p \in [1 , \infty ]$,
every $d \in \N$ and every $x = (x_1, \ldots, x_d ) \in \R^d$ the $\ell_p$-norm of $x$
is defined by 
\begin{equation*}
\norm{x} _p = 
\begin{cases}
\rbr[\big]{ \abs{x_1 } ^p + \cdots + \abs{x_d} ^p } ^{ \nicefrac{1}{p} } & \colon p < \infty \\
\max \cu{\abs{x_1}, \ldots, \abs{x_d } } & \colon p = \infty .
\end{cases}
\end{equation*}
For every $n \in \N_0 =\N \cup \cu{0} = \cu{0, 1, 2, \ldots }$ we write $\num{n} = \N \cap (0 , n ] = \cu{1, 2, \ldots, n }$.


 We now introduce the special function classes we consider throughout the remainder of this article.
Firstly, we consider for $d \in \N$ the multidimensional maximum and product functions
$m_d, p_d  \colon \R^d \to \R$ and $ \fm_d, \fp_d  \colon \R^d \to \R^d$
which
satisfy for all $x = (x_1, \ldots, x_d ) \in \R^d$ that
\begin{equation}
\label{eq:max:prod}
    \begin{split}
        m_d ( x ) & = \max \cu{x_1, \ldots, x_d} , \\
        p_d ( x ) &= \sprod_{i=1}^d x_i, \\
        \fm_d ( x ) &= ( m_1 ( x_1 ), m_2  ( x_1, x_2 ) , \ldots, m_d ( x_1,  \ldots, x_d ) ) , \\
        \fp_d ( x ) &= ( p_1 ( x_1 ) , p_2 ( x_1, x_2 ) , \ldots, p_d ( x_1, \ldots, x_d ) ).
    \end{split} 
\end{equation}

Furthermore, for every $n \in \N$,
$d_1, \ldots, d_n$,
$e_1, \ldots, e_n \in \N$,
subsets $A_i \subseteq \R^{ d_i}$,
and arbitrary functions $f_i \colon A_i \to \R^{ e_i}$,
$i \in \num{n}$,
we denote by 
$( f_1 \square f_2 \square \cdots \square f_n )
\colon A_1 \times A_2 \times \cdots \times A_n \to \R^{\sum_{i=1}^n e_i }$
the \emph{parallelization} of $f_1, \ldots, f_n$,
i.e., the function which satisfies for all $x_1 \in A_1 , \ldots, x_n \in A_n $ that
\begin{equation}
    ( f_1 \square f_2 \square \cdots \square f_n ) ( x_1, x_2, \ldots, x_n ) = (f_1 ( x_1 ) , f_2 ( x_2 ) , \ldots, f_n ( x_n ) ).
\end{equation}
Now let $Q$ be a ($d$-dimensional) \emph{hypercube}, by which we in the following mean a set of the form $Q = [a, b ] ^d \subseteq \R^d$.
Here $d \in \N$ is an arbitrary natural number which describes the input dimension and $a \in \R$, $b \in (a, \infty )$ are real numbers which describe the input domain.
For every $k\in \N$, $p \in [1 , \infty ]$,
$L \in [0 , \infty )$
we denote by $\scrP_{k, p} ( Q , L ) \subseteq \bigcup_{n \in \N} C ( Q , \R^n )$
the set of all continuous functions $f $ on $Q$
for which there exist $n \in \N$,
$d_1, \ldots, d_n \in \N$,
and functions
$f_i \colon [a,b] ^{ d_i } \to \R$, $i \in \num{n}$,
such that
\begin{equation*}
\begin{split}
	\ssuml_{i=1}^ n d_i = d,
	\qquad & \max _{ i \in \num{n}} d_i \le k, \qquad
	\forall \, i \in \num{n}, \, x, y \in [a, b ] ^{ d_i } \colon \abs{ f_i ( x ) - f_i ( y ) } \le L \norm{x - y } _p , \\
	\qqandqq 
	& f = f_1 \square f_2 \square \cdots \square f_n .
\end{split} 
\end{equation*}
That is, $f$ is the parallelization of functions defined on domains of dimension at most $k$ which are Lipschitz with respect to the $p$-norm with Lipschitz constant $L$.

We also denote
by $\bfM ( Q ) \subseteq \bigcup_{n \in \N } C ( Q , \R^n )$
the set of all continuous functions $f $
 on $ Q $
 for which there exist $n \in \N$
 and
$d_1, \ldots, d_n \in \N$
such that
$\sum_{i=1}^ n d_i = d$
and
\begin{equation*}
    f = ( m_{d_1} | _{[a , b ] ^{d_1 } } ) \square ( m_{d_2} | _ {[a, b ] ^{d_2 } } ) \square \cdots \square ( m_{d_n} | _ { [a, b ] ^{d_n } } ),
\end{equation*}
i.e., the parallelizations of the maximum functions $m_{d_i}$.

Finally, we
denote 
by $\bfP ( Q  ) \subseteq \bigcup_{n \in \N } C ( Q , \R^n )$
the set of all continuous functions $f $
 on $ Q $
 for which there exist $n \in \N$
 and
$d_1, \ldots, d_n \in \N$
such that
$\sum_{i=1}^ n d_i = d$
and
\begin{equation*}
    f = ( p_{d_1} | _{[a , b ] ^{d_1 } } ) \square ( p_{d_2} | _ {[a, b ] ^{d_2 } } ) \square \cdots \square (p_{d_n} | _ { [a, b ] ^{d_n } } ) ,
\end{equation*}
i.e., the parallelizations of the product functions $p_{d_i}$.

\subsection{Main results}
In the following we show how to approximate certain compositions of parallelized functions by DNNs without the curse of dimensionality.
Our results are somewhat similar to the previous results in Cheridito et al.~\cite{Cheridito2022_Catalog} and Beneventano et al.~\cite{BeneventanoGraeber2021}. 
While the authors of \cite{Cheridito2022_Catalog} use the general framework of catalog networks, our arguments exploit the compositional structure of the target functions more directly.
While some of our arguments rely on the results in \cite{BeneventanoGraeber2021}, a main improvement is that we consider compositions of a variable and potentially unbounded number of functions.
Specifically, in our first main result, \cref{theo:intro:1} below,
the number $\bfk ( d ) $ of functions in the composition is allowed to grow at most polynomially in the parameter $d \in \N$ which describes the dimension. 
An additional improvement in \cref{theo:intro:1} compared to \cite{BeneventanoGraeber2021} is that the functions $g_i^d$ are allowed to be parallelizations of Lipschitz functions of input dimension at most $c \in \N$ (the class $\scrP_{c, 1 }$)
instead of only $1$-dimensional Lipschitz functions.

We now present the precise statement of \cref{theo:intro:1} and, thereafter, illustrate this statement by means of several examples.

\begin{theorem} 
	\label{theo:intro:1}
	Let 
	$c \in \N$,
	for every $d \in \N$ let $\bfk ( d ) , \fd_1^d , \fd_2^d, \ldots,  \fd_{\bfk ( d ) + 1 } ^d \in \num{c d^c }$,
	for every $d \in \N$, $i \in \num{\bfk ( d ) } $
	let $Q_i^d \subseteq [ - c d^c , c d^c ] ^{ \fd_i^d }$ be a $\fd_i^d$-dimensional hypercube and
	let $g_i^d \in C ( Q_i^d , \R^{ \fd_{i + 1 } ^d  } )$ be a function,
	assume for every $d \in \N$, $i \in \num{\bfk ( d ) } $
	that
	\begin{equation} 
	\label{theo:intro:1:eq:allowed:function}
	g_i^d \in \scrP_{c, 1 } ( Q_i^d , 1 )
	\cup \bfM ( Q_i^d)
	\cup \bfP ( Q_i^d  ),
	\end{equation}
	assume for all $d \in \N$, $i \in \num{\bfk ( d ) - 1 } $
	that $g_i^d (Q_i^d ) \subseteq Q_{i+1} ^d$,
	assume for all $d \in \N$,
	$i \in \num{\bfk ( d ) }$
	with $g_i^d \in \bfP (Q_i^d )$
	that $Q_i^d  \subseteq [ - \frac{1}{8} , \frac{1}{8} ]^{ \fd_i^d}$,
	and for every $d \in \N$ let $F_d \in C( Q_1^d , \R^{ \fd _{\bfk ( d ) + 1 } ^d } )$ satisfy $F_d = g^d_{ \bfk ( d ) } \circ g^d_{\bfk ( d ) - 1 } \circ \cdots \circ g^d_{ 1 }$.
	Then there exists $K \in \N$
	such that for every $d \in \N$, $\varepsilon \in (0 , 1 ]$
	there exists $\Phi \in \bfN$
	such that
	$\cR ( \Phi ) \in C ( \R^{\fd_1^d } , \R^{\fd_{\bfk ( d ) + 1 } ^d } )$,
	\begin{equation}
	\label{theo:intro:eq:result}
	\cP ( \Phi ) \le K d^K \varepsilon^{ - 2 c }, \qqandqq
	\forall \, x \in Q_1^d \colon  \norm{ \cR ( \Phi ) ( x ) - F_d ( x ) } _ 1 
	\le \varepsilon .
	\end{equation}
\end{theorem}

The condition \cref{theo:intro:1:eq:allowed:function} in \cref{theo:intro:1}
asserts that each function $g_i^d$ in the composition is either a parallelized Lipschitz function with Lipschitz constant $1$, a parallelized maximum function, or a parallelized product function of $\fd_i^d \le c d^c$ variables. Here we think of $d \in \N$ as a parameter describing the order of magnitude of the dimensions. A particular case is $\fd_1^d = d$ for all $d \in \N$, so that the target function $F_d $ is defined on a $d$-dimensional domain.

We need to assume
for all $d \in \N$,
$i \in \num{\bfk ( d ) }$
with $g_i^d \in \bfP (Q_i^d )$
that $Q_i^d  \subseteq [ - \frac{1}{8} , \frac{1}{8} ]^{ \fd_i^d}$
 since only on sufficiently small hypercubes the parallelized product functions in $\bfP ( Q_i^d  )$ can be approximated with Lipschitz constant $1$. For Lipschitz constants greater than $1$, the composition  could potentially produce exponentially large errors and thus could not be approximated without the curse of dimensionality.
 
 The conclusion of the theorem in \cref{theo:intro:eq:result} establishes that each of the functions $F_d$ can be approximated uniformly with accuracy $\varepsilon$ on the $\fd_1^d$-dimensional hypercube $Q_1^d $ by a DNN with ReLU activation using at most $K d^K \varepsilon^{ - 2c }$ parameters. Hence the number of parameters grows at most polynomially in $d$ and the accuracy $\varepsilon$.

Next let us illustrate the statement of \cref{theo:intro:1} by some examples, which demonstrate several cases of compositions of a number of functions depending on the input dimension.

\begin{example}
A situation where \cref{theo:intro:1} applies is the family of functions 
\begin{equation*}
     [e ^{-1} , 1 ] ^d \ni x = (x_1, x_2, \ldots, x_d ) \mapsto x_1^{x_2^{ \cdots ^{x_d } } } \in \R,
     \quad d \in \N,
\end{equation*}
which can be written as $g_1^d \circ \cdots \circ g_{d-1}^d$ where 
$g_i ^d \colon [e^{-1} , 1 ] ^{i+1} \to [e^{-1} , 1 ] ^i$ defined by $g_i ^d ( x_1, \ldots, x_{i-1}, x_i, x_{i+1} ) = (x_1, \ldots, x_{i-1}, x_i^{x_{i+1} } ) $ is an element of $\scrP_{2, 1} ( [e^{-1} , 1 ] ^{i+1} , 1 )$, as can be easily verified. \cref{theo:intro:1} hence implies that these functions can be approximated by DNNs without the curse of dimensionality.
\end{example}

\begin{example}
We can also consider for an arbitrary $a \in (1, \infty )$ the family of functions
\begin{equation*}
    [1, a ] ^d \ni x = (x_1, x_2, \ldots, x_d ) \mapsto \ln ( x_1 + \ln ( x_2 + \ln ( \cdots + \ln ( x_{d } ) ) ) ) \in \R, \quad d \in \N .
\end{equation*}
Define $a_1^d , a_2^d , \ldots, a_d^d \in (1, \infty )$ by $a_d^d = a$ and $a_{i} ^d = a_{i+1} ^d + \ln ( a_{i+1} ^d )$ for $ i \in \num{d-1}$.
It is not hard to see that $a_i^d \le cd^c$ for a suitable constant $c$. Now the functions in question can be written as $g_1^d \circ \cdots \circ g_{d-1} ^d$
where $g_i^d \colon [1, a _{ i + 1 } ^d] ^{ i + 1 } \to [1 , a_{i} ^d ]^i$,
defined by $g_i^d (x_1, \ldots, x_{i-1}, x_i, x_{i+1} ) = (x_1, \ldots, x_{i-1}, x_i + \ln ( x_{i+1} ) )$,
is an element of
$\scrP_{2, 1 } ( [1, a_{ i + 1 } ^d ] ^{ i + 1 } , 1 )$.
By \cref{theo:intro:1} these functions can thus be approximated by DNNs without the curse of dimensionality.
\end{example}

\begin{example}
For arbitrary $a \in (0 , \frac{1}{8} ] $
we can consider the family of functions
\begin{multline*}
    [- a , a ] ^{d^2  } \ni x = (x_1, x_2, \ldots, x_{d^2 } )  \mapsto \\
   \rbr[\big]{\sprod_{i=1}^d x_i } \max \cu[\big]{\max\nolimits_{j=1}^d x_{d+j} ,  \rbr[\big]{ \sprod_{k=1}^d x_{2 d + k } }
   	 \max \cu{ \max\nolimits_{l=1}^d x_{3 d + l } , \cdots} } \in \R ,
    \quad d \in \N,
\end{multline*}
which is a composition of $\le d $ parallelized maximum and product functions and thus satisfies the conditions of the theorem. For example, if $d=4$ this is the function
\begin{equation*}
(x_1, \ldots, x_{16} ) \mapsto x_1 x_2 x_3 x_4 \max \cu*{ x_5, x_6, x_7, x_8, x_9 x_{10} x_{11} x_{12} \max \cu{ x_{13} , x_{14} , x_{15} , x_{16} } } .
\end{equation*}
By \cref{theo:intro:1} these functions can thus again be approximated by DNNs without the curse of dimensionality.
\end{example}

In our second main result, \cref{theo:intro:2},
the number of functions in the composition is a fixed integer $k \in \N$,
but the Lipschitz constants of the functions in the composition are allowed to depend on the dimension $d \in \N$. This is also an improvement compared to \cite{BeneventanoGraeber2021} where the maximal Lipschitz constant is a fixed number independent of $d \in \N$.
Here we allow both the parallelized maximum and product functions in $\bfM, \bfP$, as well as the generalized multidimensional maximum and product functions $\fm_d, \fp_d$.

\begin{theorem}
	\label{theo:intro:2}
	Let $c , k \in \N$,
	for every $d \in \N$ let $ \fd_1^d , \fd_2^d, \ldots,  \fd_{ k + 1 } ^d \in \num{c d^c }$,
	for every $d \in \N$, $i \in \num{ k } $
	let $Q_i^d \subseteq [ - c d^c , c d^c ] ^{ \fd_i^d }$ be a hypercube and
	let $g_i ^d \in C (Q_i^d , \R^{ \fd_{i + 1 } } )$ be a function,
	let $p \in [1 , \infty ]$,
	assume for every $d \in \N$, $i \in \num{ k } $
	that
	\begin{equation}
	\label{theo:intro:2:eq:allowed:function}
	g_i^d \in \scrP_{c, p } ( Q_i^d , c d ^c )
	\cup \cu{ \fm_{ \fd_i^d } | _ {Q_i^d } }
	\cup \cu{ \fp_{ \fd_i^d } | _ {Q_i^d } }
	\cup \bfM ( Q_i^d  )
	\cup \bfP ( Q_i^d ),
	\end{equation}
	assume for all $d \in \N$, $i \in \num{ k - 1 } $
	that $g_i^d ( Q_i^d ) \subseteq Q_{i + 1 }^d$,
	assume for all $d \in \N$,
	$i \in \num{ k }$
	with $g_i^d \in \bfP (  Q_i^d  ) \cup \cu{ \fp_{ \fd_i^d } | _ { Q_i^d } }$
	that $Q_i^d \subseteq [ - 1 , 1 ]^{ \fd_i^d }$,
	and for every $d \in \N$ let $F_d \in C( Q_1^d , \R^{ \fd _{ k + 1 } ^d } )$ satisfy $F_d = g^d_k \circ g^d_{k - 1 } \circ \cdots \circ g^d_{ 1 }$.
	Then there exists $K \in \N$
	such that for every $d \in \N$, $\varepsilon \in (0 , 1 ]$
	there exists $\Phi \in \bfN$
	such that
	$\cR ( \Phi ) \in C ( \R^{\fd_1^d } , \R^{\fd_{k + 1 } ^d } )$,
	\begin{equation}
	\cP ( \Phi ) \le K d^K \varepsilon^{ - 2 c }, 
	\qqandqq \forall \, x \in Q_1^d  \colon
	\norm{ \cR ( \Phi ) ( x ) - F_d ( x ) } _ p 
	\le \varepsilon .
	\end{equation}
\end{theorem}

Again, we need to assume for all $d \in \N$,
$i \in \num{ k }$
with
$g_i^d \in \bfP ( Q_i^d ) \cup \cu{ \fp_{ \fd_i^d } | _ { Q_i^d } }$
that
$Q_i^d \subseteq [ - 1 , 1 ] ^{ \fd_i^d }$
since on larger hypercubes the multidimensional product functions can only be approximated with a Lipschitz constant growing exponentially in the dimension (cf., e.g., \cite[Proposition 6.8]{BeneventanoGraeber2021}).

Let us also illustrate the statement of \cref{theo:intro:2} by means of a few examples.
\begin{example}
Consider the family of functions
\begin{equation*}
    [- 1 , 1 ]^d \ni x \mapsto \max \cu[\big]{(x_1) ^d , (x_1 x_2)^{d+1}, \ldots, (x_1 x_2 \cdots x_d)^{2d - 1} } \in \R,
    \quad
    d \in \N,
\end{equation*}
which can be written as a composition of three functions $g_3^d \circ g_2^d \circ g_1^d$. Here $g_1^d = \fp_d \colon [-1, 1 ]^d \to [-1, 1 ] ^d $ is the extended product function,
$g_2^d \colon [-1, 1 ] ^d \to [-1, 1 ] ^d$ is the parallelization of the component-wise functions $[-1 , 1 ] \ni y_i \mapsto (y_i)^{ d+ i - 1 } \in \R $ which have a Lipschitz constant that grows polynomially in the dimension $d$,
and $g_3^d = m_d \colon [-1, 1 ] ^d \to \R$ is the maximum function. \cref{theo:intro:2} thus implies that these functions can be approximated by DNNs without the curse of dimensionality.
\end{example}

\begin{example}
Consider for arbitrary $c \in \N$
the family of functions
\begin{equation*}
    [- c d ^c , c d ^c ] ^d \ni x \mapsto \sprod_{i=1}^d \exp \rbr*{ - i \abs{x_i } ^2 } \in \R ,
    \quad d \in \N ,
\end{equation*}
which is the composition of the product function $p_d$ and the parallelization of the one-dimensional functions $x_i \mapsto \exp \rbr*{ - i \abs{x_i } ^2 } $ which have a Lipschitz constant bounded by a polynomial in $d$ and map into the interval $[0, 1 ]$. By \cref{theo:intro:2} these functions can thus be approximated by DNNs without the curse of dimensionality.
\end{example}

\begin{example}
	Consider for arbitrary $c \in \N$
	the family of functions
	\begin{equation*}
	[ -c d ^c , c d ^c ] ^{ 3 d } \ni x \mapsto \max_{l=1}^d \cos ( l x _{ 3 l - 2 } + l ^2 x _ { 3 l - 1 } + l ^3  x_{ 3 l } ) \in \R , 
	\quad d \in \N ,
	\end{equation*}
	which can again be approximated by DNNs without the curse of dimensionality by \cref{theo:intro:2}, since the functions $\R^3 \ni ( y_1, y_2, y_3 ) \mapsto \cos ( l y_1 + l ^2 y_2 + l^3 y_3 ) \in \R $, $1 \le l \le d$, have a Lipschitz constant bounded by a polynomial in $d$.
\end{example}

The above examples illustrate the type of compositional functions that can now be approximated by DNNs with number of parameters growing at most polynomially in the input dimension and the approximation accuracy.

The remainder of this article is organized as follows. In \cref{section:dnn:approx} we establish abstract DNN approximation results for certain function compositions. 
Afterwards, in \cref{section:approx:function:classes} we prove that under suitable assumptions the parallelized Lipschitz, product, and maximum functions introduced above can be approximated without the curse of dimensionality. We then combine this with the abstract results from \cref{section:dnn:approx} to establish the main theorems.

\section{DNN approximation of composite functions}
\label{section:dnn:approx}

Throughout this section we employ the notation introduced in \cref{subsec:dnn:notation}.
The main results of this section are the two abstract approximation results for compositions of functions in \cref{prop:main:abstract:comp1,prop:main:abstract:comp2} which cover the cases of a polynomially growing number of functions with Lipschitz constant $1$ and of a fixed number of functions with polynomially growing Lipschitz constants, respectively.
To prove these we require some technical preparations.

\subsection{DNN approximation cost}
We first define in \cref{def:cost} our version of the approximation cost of continuous functions by DNNs and then establish some basic properties of the approximation cost.
\cref{def:cost} is inspired by Beneventano et al.~\cite[Definition 3.2]{BeneventanoGraeber2021}, but is in a sense more general since it allows us to consider different $\ell_p$-norms on $\R^n$.

\begin{definition}[Cost of DNN approximations]
\label{def:cost}
For every
$p \in [1 , \infty ]$,
$m, n \in \N$,
$D \subseteq \R^m$,
$f \in C ( D , \R^n )$,
$L , \varepsilon \in [0 , \infty )$
we denote by $\cost{p} ( f , L , \varepsilon ) \in \N \cup \cu{\infty}$
the extended real number given by
\begin{equation} 
\label{eq:def:cost}
\begin{split} 
   \cost{p} ( f , L , \varepsilon ) & =
   \min \Big( \Big\{ N \in \N \colon \exists \, \Phi \in \bfN \colon \big[ ( \cR ( \Phi ) \in C ( \R^m , \R^n ) ) \\
  & \wedge ( N = \cP ( \Phi ) ) 
    \wedge \rbr[\big]{ \sup\nolimits_{x \in D } \norm{ \cR ( \Phi ) ( x ) - f ( x ) }_p \le \varepsilon } \\
   & \wedge \rbr[\big]{ \forall \, x,y \in D \colon \norm{ \cR ( \Phi ) ( x ) - \cR ( \Phi ) ( y ) }_p \le L \norm{ x - y }_p } \big]
   \Big\} \cup \cu{\infty} \Big) .
   \end{split}
\end{equation}
\end{definition}

That is, $ \cost{p} ( f , L , \varepsilon ) $ is the minimal number of parameters of a DNN which can approximate $f$ up to accuracy $\varepsilon$ in the $\ell_p$-norm while being $L$-Lipschitz on the domain of $f$ with respect to the $\ell_p$-norm. If such a DNN does not exist then the cost is defined to be infinite.
Controlling the Lipschitz constant will be important when estimating the propagation of the approximation error through compositions; see \cref{prop:approx:composition} below.
 
We next restate in \cref{lem:cost:mono} below the monotonicity of the approximation cost with respect to the parameters $L, \varepsilon$ and the domain $D$, which was established, e.g., in \cite[Lemma 3.8]{BeneventanoGraeber2021} for the case $p=2$. The case of general $p$ is entirely analogous.

\begin{lemma}
\label{lem:cost:mono}
Let $p \in [1 , \infty ]$,
$m, n \in \N$,
$D \subseteq \R^m$,
$E \subseteq D$,
$f \in C ( D , \R^n )$,
$L_1, L_2 , \varepsilon_1, \varepsilon_2 \in [0 , \infty )$
satisfy $L_1 \le L_2$ 
and $\varepsilon_1 \le \varepsilon _2$.
Then
\begin{equation}
    \cost{p} ( f, L_1, \varepsilon_1 ) \ge \cost{p} ( f | _E , L_2 , \varepsilon _ 2 ).
\end{equation}
\end{lemma}

In the next lemma we show for arbitrary $p , q \in [1 , \infty ]$
how to estimate the cost with respect to the $\ell_q$-norm against the cost with respect to the $\ell_p$-norm.
For this we need to rescale the Lipschitz constant $L$ and the approximation accuracy $\varepsilon$ by factors depending on the input and output dimensions $m, n$.

\begin{lemma}
[Cost with respect to different norms]
\label{cost:p:norm}
Let
$p, q \in [1 , \infty ]$,
$m, n \in \N$,
$D \subseteq \R^m$,
$f \in C ( D , \R^n )$,
$L , \varepsilon \in [0 , \infty )$.
Then
\begin{equation}
    \cost{q} \rbr[\big]{ f , \max \cu[\big]{m ^{ \nicefrac{1}{p} - \nicefrac{1}{q} } , n ^{ \nicefrac{1}{q} - \nicefrac{1}{p} } } L , \max \cu[\big]{ n ^{ \nicefrac{1}{q} - \nicefrac{1}{p} } , 1 } \varepsilon }
     \le \cost{p} ( f , L , \varepsilon ).
\end{equation}
\end{lemma}

In \cref{cost:p:norm} above we use the convention that $ \frac{1}{\infty} = 0$.

\begin{cproof}{cost:p:norm}
	Throughout this proof we use the well-known fact that for all $k \in \N$, $x \in \R^k $, $1 \le s \le t \le \infty $ it holds that $\norm{x} _ t \le \norm{x} _s \le k^{ \nicefrac{1}{s} - \nicefrac{1}{t} } \norm{x} _t$.
Assume w.l.o.g.~that 
$\cost{p} ( f , L , \varepsilon ) < \infty $.
\cref{eq:def:cost} therefore
assures that there exists $\Phi \in \bfN$
which satisfies 
\begin{multline*}
	\cR ( \Phi ) \in C ( \R^m , \R^n ), \qquad
	\forall \, x \in D \colon \norm{\cR ( \Phi ) ( x ) - f ( x ) } _p \le \varepsilon ,
	 \\
	\forall \, x,y \in D \colon \norm{ \cR ( \Phi ) ( x ) - \cR ( \Phi ) ( y ) }_p \le L \norm{ x - y }_p,
	\qqandqq 
	\cP ( \Phi ) = \cost{p} ( f , L , \varepsilon ) .
\end{multline*}
This yields for all $x \in D$ that 
\begin{equation*}
    \norm{\cR ( \Phi ) ( x ) - f ( x ) } _ q 
    \le \max \cu[\big]{ n ^{ \nicefrac{1}{q} - \nicefrac{1}{p} } , 1 } \norm{ \cR ( \Phi ) ( x ) - f ( x ) }_p 
    \le \max \cu[\big]{ n ^{ \nicefrac{1}{q} - \nicefrac{1}{p} } , 1 } \varepsilon.
\end{equation*}
Furthermore, if $q \ge p$ we obtain for all $x, y \in D$ that
\begin{equation*}
    \begin{split}
        \norm{ \cR ( \Phi ) ( x ) - \cR ( \Phi ) ( y ) }_q
        & \le \norm{ \cR ( \Phi ) ( x ) - \cR ( \Phi ) ( y ) }_p
        \le L \norm{x - y }_p \\
        & \le m ^{ \nicefrac{1}{p} - \nicefrac{1}{q} } L \norm{x - y } _ q.
    \end{split}
\end{equation*}
On the other hand, if $q \le p$ we get for all $x, y \in D$ that
\begin{equation*}
    \begin{split}
        \norm{ \cR ( \Phi ) ( x ) - \cR ( \Phi ) ( y ) }_q
        & \le n ^{ \nicefrac{1}{q} - \nicefrac{1}{p} } \norm{ \cR ( \Phi ) ( x ) - \cR ( \Phi ) ( y ) }_p
        \le n ^{ \nicefrac{1}{q} - \nicefrac{1}{p} } L  \norm{x - y }_p \\
        &\le n ^{ \nicefrac{1}{q} - \nicefrac{1}{p} } L \norm{x - y } _ q.
    \end{split}
\end{equation*}
In any case, we have $\forall \, x,y \in D \colon \norm{ \cR ( \Phi ) ( x ) - \cR ( \Phi ) ( y ) }_q
\le \max \cu{ m ^{ \nicefrac{1}{p} - \nicefrac{1}{q} } , n ^{ \nicefrac{1}{q} - \nicefrac{1}{p} }  } L \norm{x - y } _ q $.
\end{cproof}

\begin{remark}
	Note that a particularly simple case of \cref{cost:p:norm} arises if $n=1$, i.e. the output dimension of the considered functions is $1$, since then all $\ell_p$-norms on the output space agree. In this case, we obtain 
	\begin{equation*}
	\cost{q} \rbr[\big]{ f , \max \cu[\big]{ m ^{ \nicefrac{1}{p} - \nicefrac{1}{q} } , 1 } L ,  \varepsilon } \le \cost{p} ( f , L , \varepsilon ).
	\end{equation*}
	In particular, if $p \ge q$ we simply get $\cost{q} ( f , L , \varepsilon ) \le \cost{p} ( f , L , \varepsilon )$.
\end{remark}

We next establish in \cref{lem:cost:clipped} an auxiliary result which shows that if the values of a function are contained in a hypercube $Q$, the output values of the approximating DNN can be clipped to be contained in $Q$ as well with only a moderate increase in the number of parameters. This result will be useful when composing multiple DNNs.
\cref{lem:cost:clipped} is a generalization of \cite[Lemma 3.7]{BeneventanoGraeber2021},
which covers the case $p=2$ and $Q = [-R , R ] ^n$.
The proof is very similar.

\begin{lemma} \label{lem:cost:clipped}
	Let
	$p \in [1 , \infty ]$,
	$m, n \in \N$,
	$D \subseteq \R^m$,
	$f \in C ( D , \R^n )$,
	$L , \varepsilon \in [0 , \infty )$ satisfy 
	$\cost{p} ( f , L , \varepsilon) < \infty $,
	let $Q \subseteq \R^n$ be a hypercube, and assume $f ( D ) \subseteq Q$.
	Then there exists $\Phi \in \bfN$
	which satisfies
	\begin{multline}
	\cR ( \Phi ) \in C ( \R^m , \R^n ), \qquad
	\forall \, x \in D \colon \cR ( \Phi ) ( x ) \in Q , \\
	\forall \, x \in D \colon \norm{ \cR ( \Phi ) ( x ) - f ( x ) } _p \le \varepsilon, \qquad
	\forall \, x , y \in D \colon \norm{\cR ( \Phi ) ( x ) - \cR ( \Phi ) ( y ) } _ p \le L \norm{x - y } _p , \\
	\qqandqq
	\cP ( \Phi ) \le \cost{p} ( f , L , \varepsilon) + 2 n ( n + 1 ) .
	\end{multline}
\end{lemma}

\begin{cproof}{lem:cost:clipped}
	\Nobs that the assumption that $\cost{p} ( f , L , \varepsilon) < \infty $
	ensures that there exists $\Psi \in \bfN$
	which satisfies
	\begin{multline}
	\label{lem:cost:clipped:eq1}
	\cR ( \Psi ) \in C ( \R^m , \R^n ),  \qquad
	\forall \, x \in D \colon \norm{ \cR ( \Psi ) ( x ) - f ( x ) } _p \le \varepsilon, \\
	\forall \, x , y \in D \colon \norm{\cR ( \Psi ) ( x ) - \cR ( \Psi ) ( y ) } _ p \le L \norm{x - y } _p ,
	\qqandqq
	\cP ( \Psi ) = \cost{p} ( f , L , \varepsilon) .
	\end{multline}
	Next assume $Q = [ a , b ] ^n$
	and let $\fC \colon \R^n \to \R^n$ denote the multidimensional clipping function which satisfies for all
	$x = (x_1, \ldots , x_n ) \in \R^n$
	that
	\begin{equation*}
	\fC ( x ) = ( \max \cu{ a , \min \cu{x _1 , b } } , \ldots, \max \cu{ a , \min \cu{x _n , b } } ) .
	\end{equation*}
	\Nobs that for all $x , y \in \R^n$ we have that $\fC ( x ) \in Q$
	and $\norm{ \fC ( x ) - \fC ( y ) } _p \le \norm{x - y } _ p$.
	Furthermore, \cite[Lemma 2.32]{BeneventanoGraeber2021}
	ensures that there exists a network $\Theta \in \bfN$
	with architecture $(\ell_0, \ell_1, \ell_2, \ell_3) = (n, n, n, n)$
	which satisfies for all $x \in \R^n$ that $\cR ( \Theta ) ( x ) = \fC ( x ) $.
	Now define $\Phi \in \bfN$ as the composition $\Phi = \Theta \circ \Psi $ (see \cite[Definition 2.15]{BeckJentzenKuckuck2022} for the formal definition).
	Since $\cO (\Psi ) = n$, it is not hard to see that $\cP ( \Phi )= \cP ( \Psi ) + 2 n ( n + 1 )$.
	Furthermore, we have for all $x \in D$ that
	$\cR ( \Phi ) ( x ) = ( \cR ( \Theta ) \circ \cR ( \Psi ) ) ( x ) = \fC ( \cR ( \Psi ) ( x ) ) \in Q$.
	Moreover, since $f(D) \subseteq Q$ it follows for all $x \in D$ that $\fC ( f ( x ) ) = f ( x )$ and therefore
	\begin{equation*}
	\norm{\cR ( \Phi ) ( x ) - f ( x ) }_p = \norm{\fC ( \cR ( \Psi ) ( x ) ) - \fC ( f ( x ) ) }_p
	\le \norm{\cR ( \Psi ) ( x ) - f ( x ) } _ p \le \varepsilon .
	\end{equation*}
	In the same way, for all $x, y \in D$ the Lipschitz property of $\fC$ implies
	\begin{equation*}
	\begin{split}
	\norm{\cR ( \Phi ) ( x ) - \cR ( \Phi ) ( y ) } _ p 
	& = \norm{ \fC ( \cR ( \Psi ) ( x ) ) - \fC ( \cR ( \Psi ) ( y ) ) }_p
	\le \norm{\cR ( \Psi ) ( x ) - \cR ( \Psi ) ( y ) } _p \\
	& \le L \norm{x - y } _p.
	\end{split}
	\end{equation*}
	This shows that $\Phi$ satisfies \cref{lem:cost:clipped:eq1}.
\end{cproof}

\subsection{DNN approximations of compositions}

We now establish in \cref{prop:approx:composition} an upper bound for the approximation cost of the composition of an arbitrary number $n \in \N$ of functions. It is a generalization of \cite[Proposition 3.18]{BeneventanoGraeber2021}, which is essentially the case $n=2$ (note that the case $n=1$ is trivial).
The proof is similar, using \cite[Lemma 6.5]{BeneventanoGraeber2021} to estimate the error induced by the composition.

\begin{prop}
	\label{prop:approx:composition}
	Let 
	$p \in [1 , \infty ]$,
	$n \in \N$,
	let $d_1, d_2, \ldots, d_{n+1} \in \N$,
	$ L_1, L_2, \ldots, L_n \in ( 0 , \infty )$,
	for every $i \in \num{n}$ let
	$Q_i \subseteq \R^{d_i}$ be a $d_i$-dimensional hypercube,
	and let $f_i \in C ( Q_i , \R^{ d_{i+1}} ) $, $i \in \num{ n }$,
	satisfy for all $i \in \num{n - 1 }$
	that
	$f_i ( Q_i ) \subseteq Q_{i+1} $.
	Then it holds for all $\varepsilon \in [ 0 , \infty )$
	that
	\begin{equation}
	\begin{split}
	& \cost{p} \rbr[\big]{ f_n \circ f_{n-1} \circ \cdots \circ f_1 , \sprod_{i=1}^n L_i , \varepsilon }
	\\
	& \le 6 \ssuml_{i=2}^n d_i ( d_i + 1 ) 
	+ 3 \ssuml_{i=1}^n \cost{p} \rbr[\big]{ f_i, L_i, \varepsilon n^{-1}  ( \sprod_{j = i + 1 } ^n L_j ) ^{-1} } .
	\end{split}
	\end{equation}
\end{prop}

\begin{cproof}{prop:approx:composition}
Throughout this proof let
$\varepsilon \in [ 0 , \infty )$ 
and assume w.l.o.g.~that $n \ge 2$ and
\begin{equation}
\label{prop:approx:composition:eq1}
    \forall \, i \in \num{n} \colon
    \cost{p} \rbr[\big]{ f_i, L_i, \varepsilon n^{-1}  ( \sprod_{j = i + 1 } ^n L_j ) ^{-1} } < \infty .
\end{equation}
Combining \cref{prop:approx:composition:eq1} with \cref{lem:cost:clipped}
 implies that there exist
$\Phi_i \in \bfN$, $i \in \num{n}$,
which satisfy for all $i \in \num{ n }$ that
\begin{multline}
\label{prop:approx:composition:eq2}
    \cR ( \Phi_i ) \in C ( \R^{ d_i}, \R^{ d_{i + 1 } } ),
    \qquad
    \forall \, x \in Q_i \colon \cR ( \Phi_i ) ( x ) \in Q_{i+1} , \\
    \forall \, x \in Q_i \colon \norm{ \cR ( \Phi_i ) ( x ) - f_i ( x ) } _p \le \varepsilon n^{-1}  ( \sprod_{j = i + 1 } ^n L_j ) ^{-1} \norm{x - y }_p , \\
    \forall \, x,y \in Q_i \colon \norm{\cR ( \Phi_i )( x ) - \cR ( \Phi_i ) ( y ) } _p \le L_i \norm{x - y } _p , \\
    \qqandqq \cP ( \Phi_i ) = \cost{p} \rbr[\big]{ f_i, L_i, \varepsilon n^{-1}  ( \sprod_{j = i + 1 } ^n L_j ) ^{-1} } + 2 d_{i + 1 } ( d_{i + 1 } + 1 ) .
\end{multline}
Now define $\Phi \in \bfN$
as the composition $\Phi = \Phi_n \circ \bfI_{d_{n}} \circ \Phi_{n-1} \cdots \circ \bfI_{d_2} \circ \Phi_1$.
Here $\bfI_{d_i} \in \bfN$ denotes a suitable identity network for $\R^ { d_ i}$, i.e., it satisfies $\forall \, x \in \R^{ d_i }  \colon \cR ( \bfI_{d_i } ) ( x ) = x$ (cf., e.g., \cite[Definition 2.12]{BeneventanoGraeber2021}).
\Nobs that \cite[Proposition 2.19]{BeneventanoGraeber2021} and \cref{prop:approx:composition:eq2} imply that
\begin{equation*}
\begin{split}
    \cP ( \Phi ) \le 3 \ssuml_{i=1}^n \cP ( \Phi_i )
     \le 6 \ssuml_{i=2}^n d_i ( d_i + 1 ) 
        & + 3 \ssuml_{i=1}^n \cost{p} \rbr[\big]{ f_i, L_i, \varepsilon n^{-1}  ( \sprod_{j = i + 1 } ^n L_j ) ^{-1} } .
\end{split}
\end{equation*}
Furthermore, we have $\cR ( \Phi ) = \cR ( \Phi_n ) \circ \cdots \circ \cR ( \Phi_1 )$, and thus $\cR ( \Phi ) \in C ( \R^{ d_1 } , \R^{ d_{n + 1 } } )$ is Lipschitz continuous with Lipschitz constant $\prod_{i=1}^n L_i $.
Finally, \cite[Lemma 6.5]{BeneventanoGraeber2021} ensures that 
\begin{equation*}
    \begin{split}
       & \sup\nolimits_{x \in Q_1 } \norm{ \cR ( \Phi ) ( x ) - ( f_n \circ f_{n-1} \circ \cdots \circ f_1 ) ( x ) }_p \\
       & \le \ssuml_{i=1}^n \br*{ ( \sprod_{j = i + 1 } ^n L_j ) \varepsilon n^{-1}  ( \sprod_{j = i + 1 } ^n L_j ) ^{-1}  } =  \ssuml_{i=1}^n ( \varepsilon n^{-1} ) = \varepsilon .
    \end{split}
\end{equation*}
\end{cproof}

As a consequence, we obtain in \cref{cor:approx:composition} a general result regarding approximations of compositions of functions without the curse of dimensionality. 

\begin{cor}
\label{cor:approx:composition}
Let
$p \in [1 , \infty ]$,
$n \in \N$,
let $d_1, d_2, \ldots, d_{n+1} \in \N$,
$L_1, L_2, \ldots, L_n \in ( 0 , \infty )$,
for every $i \in \num{n}$ let
$Q_i \subseteq \R^{ d_i }$ be a hypercube and
let $f_i \in C( Q_i , \R^{ d_{i+1}} ) $ be a function,
assume for all $i \in \num{n - 1}$
that
$f_i ( Q_i ) \subseteq Q_{i+1} $,
let $\fd = \max \cu{d_1 , \ldots , d_n }$,
$\fL = \max \cu{L_1 , \ldots, L_n , 1 }$,
and assume that there exists $c \in \N$ which satisfies for all $i \in \num{n}$, $\delta \in ( 0 , 1 ] $
that 
$\cost{p} ( f_i, L_i, \delta ) \le c d_i^c L_i^c \delta^{-c}$.
Then it holds for all $\varepsilon \in (0 , 1 ] $ that
\begin{equation}
     \begin{split}
        \cost{p} \rbr[\big]{ f_n \circ f_{n-1} \circ \cdots \circ f_1 , \sprod_{i=1}^n L_i , \varepsilon }
        & \le 12 n \fd ^2 + 3 c n ^{ c + 1 } \fd ^c \fL ^{c n } \varepsilon ^{ - c } .
    \end{split}
\end{equation}
\end{cor}

\begin{cproof}{cor:approx:composition}
Applying \cref{prop:approx:composition} and the monotonicity from \cref{lem:cost:mono}
yields for all $\varepsilon \in ( 0 , 1 ]$ that
\begin{equation*}
    \begin{split}
         & \cost{p} \rbr[\big]{ f_n \circ f_{n-1} \circ \cdots \circ f_1 , \sprod_{i=1}^n L_i , \varepsilon } \\
         	& \le 6 \ssuml_{i=2}^n d_i ( d_i + 1 ) 
         + 3 \ssuml_{i=1}^n \cost{p} \rbr[\big]{ f_i, L_i, \varepsilon n^{-1}  ( \sprod_{j = i + 1 } ^n L_j ) ^{-1} } \\
         & \le 6 \ssuml_{i=1}^n 2 \fd ^2 + 3 \ssuml_{i=1}^n \cost{p} \rbr[\big]{ f_i, L_i, \varepsilon n^{-1}  \fL ^{n - i } } 
         \le
          12 n \fd ^2 + 3 \ssuml_{i=1}^n c d_i^c L_i^c \varepsilon^{-c}  n^{c}  \fL ^{ c ( n - i ) } \\
         & \le 
          12 n \fd ^2 + 3 \ssuml_{i=1}^n c \fd ^c n ^c \varepsilon^{-c} \fL ^{ c ( n - i + 1 )  } 
         \le 12 n \fd ^2 + 3 c \fd^c n^{ c + 1 } \varepsilon ^{ - c } \fL ^{ c n } .
    \end{split}
\end{equation*}
\end{cproof}

\begin{remark}
\label{remark:cost:cod}
In \cref{cor:approx:composition} we want the right-hand side to grow at most polynomially in the dimension $\fd$.
In the case $\fL = 1$, i.e.~$\max \cu{L_1, \ldots, L_n} \le 1$,
the number $n$ of functions in the composition is allowed to grow polynomially in $\fd$,
i.e., we can have $n \le c \fd ^c$ for some constant $c \in \N$.
If the maximal Lipschitz constant $\fL$ is larger than $1$ and the number $n$ of functions in the composition is a fixed constant we will also obtain an upper bound that avoids the curse of dimensionality.
\end{remark}

\subsection{Abstract approximation results without the curse of dimensionality}

In this subsection we establish two abstract approximation results for composite functions from which \cref{theo:intro:1,theo:intro:2} in the introduction will follow.
Both \cref{prop:main:abstract:comp1} and \cref{prop:main:abstract:comp2} are simple consequences of \cref{cor:approx:composition}.
First, in \cref{prop:main:abstract:comp1} we consider a composition of $\cO ( d^c )$ functions, where $d \in \N$ is a dimensionality parameter, each of which can be approximated with Lipschitz constant $1$ without the curse of dimensionality.

\begin{prop}
	\label{prop:main:abstract:comp1}
	Let $c \in \N$,
	$p \in [1 , \infty ]$,
	for every $d \in \N$ let $ \fd_1^d , \fd_2^d, \ldots,  \fd_{\bfk ( d ) + 1 } ^d \in \num{d ^c }$ and $\bfk ( d ) \in \num { c d^c }$,
	for every $d \in \N$, $i \in \num{\bfk ( d ) } $
	let $Q_i^d \subseteq \R^{\fd ^ d_i }$ be a $\fd ^d_i$-dimensional hypercube and
	let $g_i ^d \in C (Q_i^d , \R^{ \fd^d_{i + 1 } } )$ be a function,
	assume for all $d \in \N$, $i \in \num{\bfk ( d ) } $, $\varepsilon \in (0 , 1 ]$
	that
	\begin{equation*}
	\cost{p} ( g_i^d , 1 , \varepsilon ) \le c d^c \varepsilon ^{ - c } ,
	\end{equation*}
	assume for all $d \in \N$, $i \in \num{\bfk ( d ) - 1 } $ that
	$g_i^d ( Q_i^d ) \subseteq Q_{i+1} ^d  $,
	and for every $d \in \N$ let $F_d \in C( Q_1^d  , \R^{ \fd _{\bfk ( d ) + 1 } ^d } )$ satisfy $F_d = g^d_{ \bfk ( d ) } \circ g^d_{ \bfk ( d ) - 1 } \circ \cdots \circ g^d_{ 1 }$.
	Then there exists $K \in \N$
	such that for every $d \in \N$, $\varepsilon \in (0 , 1 ]$
	there exists $\Phi \in \bfN$
	such that
	$\cR ( \Phi ) \in C ( \R^{\fd_1^d } , \R^{\fd_{ \bfk ( d ) + 1 } ^d } )$,
	\begin{equation}
	\cP ( \Phi ) \le K d^K \varepsilon^{ - c }, \qqandqq
	\forall \, x \in Q_1^d   \colon \norm{ \cR ( \Phi ) ( x ) - F_d ( x ) } _p  \le \varepsilon .
	\end{equation}
\end{prop}

\begin{cproof}{prop:main:abstract:comp1}
For every $d \in \N$,
applying \cref{cor:approx:composition} (with $L_i \with 1$, $n \with \bfk ( d )$, $d_i \with \fd_i^d $)
yields for all
$\varepsilon \in (0 , 1 ]$
that
\begin{equation*}
    \begin{split}
        \cost{p} ( F_d , 1 , \varepsilon )
        &= \cost{p} ( g^d _ { \bfk ( d ) } \circ g ^d _ { \bfk ( d ) - 1 } \circ \cdots \circ g_1 ^d , 1 , \varepsilon ) \\
        & \le 12 \bfk ( d ) d ^{ 2 c } + 3 c ( \bfk ( d ) ) ^{ c + 1 } d ^{ c ^2 } \varepsilon ^{ - c } \\
        & \le 12 c d ^{ 3  c } + 3 c^ { c + 2 } d ^{ c ( 2 c + 1 ) } \varepsilon ^{ - c }
        \le 15 c ^{ c + 2 } d ^{ c ( 2 c + 1 ) } \varepsilon ^{ - c } .
    \end{split}
\end{equation*}
\end{cproof}

Next, in \cref{prop:main:abstract:comp2} we consider compositions of a constant number $k \in \N$ of functions, each of which can be approximated with Lipschitz constant $\cO ( d^c )$ without the curse of dimensionality.

\begin{prop}
	\label{prop:main:abstract:comp2}
	Let $c , k \in \N$,
	$p \in [1 , \infty ]$,
	for every $d \in \N$ let $ \fd_1^d , \fd_2^d, \ldots,  \fd_{ k + 1 } ^d \in \num{d ^c }$,
	for every $d \in \N$, $i \in \num{ k } $
	let $Q_i^d \subseteq \R^{\fd ^ d_i }$ be a $\fd ^d_i$-dimensional hypercube and
	let $g_i ^d \in C (Q_i^d , \R^{ \fd^d_{i + 1 } } )$ be a function,
	assume for all $d \in \N$, $i \in \num{ k } $, $\varepsilon \in (0 , 1 ]$
	that
	\begin{equation*}
	\cost{p} ( g_i^d , c d ^c , \varepsilon ) \le c d^c \varepsilon ^{ - c } ,
	\end{equation*}
	assume for all $d \in \N$, $i \in \num{ k  - 1 } $ that
	$g_i^d ( Q_i^d ) \subseteq Q_{i+1} ^d $,
	and for every $d \in \N$ let $F_d \in C( Q_1^d , \R^{ \fd _{ k + 1 } ^d } )$ satisfy $F_d = g^d_k \circ g^d_{ k - 1 } \circ \cdots \circ g^d_{ 1 }$.
	Then there exists $K \in \N$
	such that for every $d \in \N$, $\varepsilon \in (0 , 1 ]$
	there exists $\Phi \in \bfN$
	such that $\cR ( \Phi ) \in C ( \R^{\fd_1^d } , \R^{\fd_{k + 1 } ^d } )$,
	\begin{equation}
	\cP ( \Phi ) \le K d^K \varepsilon^{ - c }, \qqandqq 
	\forall \, x \in Q_1^d  \colon \norm{ \cR ( \Phi ) ( x ) - F_d ( x ) } _ p  \le \varepsilon .
	\end{equation}
\end{prop}

\begin{cproof}{prop:main:abstract:comp2}
For every $d \in \N$,
applying \cref{cor:approx:composition} (with $L_i \with c d^c$, $n \with  k $, $d_i \with \fd_i^d $)
shows for all $\varepsilon \in (0 , 1 ]$
that
\begin{equation*}
    \begin{split}
        \cost{p} ( F_d , ( c d ^c ) ^k , \varepsilon )
        & = \cost{p} ( g^d_k \circ g^d _{k-1} \circ \cdots \circ g^d_1  , (c d ^c ) ^k , \varepsilon ) \\
        &\le 12 k d ^{ 2 c } + 3 c k ^{ c + 1 } d ^{ c ^2 } (c d ^c ) ^{ c k } \varepsilon ^{ - c} \\
        & = 12 k d ^{ 2 c } + 3 c^{ c k + 1 } k ^{ c + 1 } d ^{ c ^2 ( 1 + k ) } \varepsilon ^{ - c }
        \le 15 k ^{ c + 1 } c ^{ c k + 1 } d ^{ c ^2 ( 1 + k ) } \varepsilon ^{ - c } .
    \end{split}
\end{equation*}
\end{cproof}

\section{DNN approximation of specific function classes}
\label{section:approx:function:classes}

In this section
we establish approximation results for the specific functions introduced in \cref{subsec:dnn:notation},
which allows us to apply \cref{prop:main:abstract:comp1,prop:main:abstract:comp2} to these functions.
At the end of this section we employ our approximation results to prove \cref{theo:intro:1,theo:intro:2} from the introduction.

\subsection{DNN approximations of parallelizations}

The first step is to show in \cref{lem:approx:parallel} how parallelizations of the form $f_1 \square f_2 \square \cdots \square f_n$, as defined in \cref{subsec:dnn:notation}, can be approximated efficiently by DNNs.
For this we employ the parallelized DNN architecture from Cheridito et al.~\cite{Cheridito2022_Catalog}, which allows one to approximate parallelizations efficiently.
Specifically, we use \cref{lem:parallel:catalog} below,
which is a reformulation of \cite[Proposition 5]{Cheridito2022_Catalog} in the ReLU case, where $c=2$ (in the notation of \cite[Proposition 5]{Cheridito2022_Catalog}).
\begin{lemma}
\label{lem:parallel:catalog}
Let $n \in \N$,
$\Phi_1, \ldots, \Phi_n \in \bfN$.
Then there exists $\Phi \in \bfN$
which satisfies $\cR ( \Phi ) = \cR ( \Phi_1 ) \square \cdots \square \cR ( \Phi_n )$
and
\begin{equation}
    \cP ( \Phi ) \le \tfrac{11}{4} n ^2 \rbr*{ \max \nolimits_{i \in \num{n} } \max \cu*{ \cI ( \Phi_i ) , \cO ( \Phi_i ) } } ^2 \ssuml_{i=1}^n \cP ( \Phi_i ) .
\end{equation}
\end{lemma}

\begin{lemma}
\label{lem:approx:parallel}
Let 
$p \in [ 1 , \infty ]$,
$n \in \N$,
$a \in \R$,
$b \in ( a , \infty )$,
$L , \varepsilon \in [0 , \infty )$,
$d_1, \ldots, d_n, \allowbreak
e_1, \ldots, e_n \in \N$,
and for every $i \in \num{n}$ let
$f_i \in C ( [a , b ] ^{ d_i } , \R ^{ e_i } )$.
Then
\begin{equation} 
    \cost{p} \rbr*{  f_1 \square f_2 \square \cdots \square f_n , L , \varepsilon }
    \le \tfrac{11}{4} n ^2 \rbr*{ \max\nolimits_{i \in \num{n} } \max \cu{d_i, e_i } } ^2 \ssuml_{i=1}^n \cost{p} \rbr*{ f_i , L , n ^{ - \nicefrac{1}{p} } \varepsilon } .
\end{equation}
\end{lemma}

\begin{cproof2}{lem:approx:parallel}
Throughout this proof assume without loss of generality that 
\begin{equation}
\label{lem:approx:parallel:eq1}
    \forall \, i \in \num{n} \colon \cost{p} \rbr[\big]{ f_i , L , n ^{ - \nicefrac{1}{p} } \varepsilon }  < \infty .
\end{equation}
\Nobs that \cref{lem:approx:parallel:eq1}
and \cref{eq:def:cost}
imply that there exist $\Phi_i \in \bfN$, $i \in \num{n}$,
which satisfy for all $i \in \num{n}$
that
\begin{multline*}
    \cR ( \Phi_i ) \in C ( \R^{ d_i}, \R^{ e_i } ),
    \qquad
    \forall \, x \in [a , b ] ^{d_i } \colon \norm{ \cR ( \Phi_i ) ( x ) - f_i ( x ) } _p \le \varepsilon n^{- \nicefrac{1}{p} }  , \\
    \forall \, x,y \in [a, b]^{ d_i} \colon \norm{\cR ( \Phi_i )( x ) - \cR ( \Phi_i ) ( y ) } _p \le L \norm{x - y } _p , \\ 
    \text{and} \qquad 
    \cP ( \Phi_i ) = \cost{p} \rbr[\big]{ f_i, L, \varepsilon n^{- \nicefrac{1}{p} } } .
\end{multline*}
Let $\Phi \in \bfN$ be the parallelization of $\Phi_1, \ldots, \Phi_n$ given by \cref{lem:parallel:catalog}. That is,
we have 
$\cR ( \Phi ) = \cR ( \Phi_1 ) \square \cdots \square \cR ( \Phi_n ) \in C ( \R ^{ \sum_{i=1}^n d_i } , \R ^{ \sum_{i=1}^n e_i } )$
and
\begin{equation*}
\begin{split}
    \cP ( \Phi ) 
    &\le \tfrac{11}{4} n ^2 \rbr*{ \max\nolimits_{i \in \num{n} }  \max \cu{d_i, e_i } } ^2 \ssuml_{i=1}^n \cP ( \Phi_i ) \\
    & = \tfrac{11}{4} n ^2 \rbr*{ \max\nolimits_{i \in \num{n} } \max \cu{d_i, e_i } } ^2 \ssuml_{i=1}^n \cost{p} ( f_i , L , n ^{ - \nicefrac{1}{p} } \varepsilon ) .
\end{split}
\end{equation*}
We now consider two cases depending on whether $p \in [1 , \infty )$ or $p = \infty $.
\begin{case}
Assume $p < \infty $.
Thus we get for all $x = (x_1, \ldots, x_n) \in [a, b ] ^{ \sum_{i=1}^n d_i } $ with $ \forall \, i \in \num{n} \colon x_i \in [a, b ] ^{ d_i }$
that
\begin{equation*}
    \begin{split}
        \norm{\cR ( \Phi ) ( x ) - ( f_1 \square f_2 \square \cdots \square f_n ) ( x ) } _ p ^p
        & = \ssuml_{i=1}^n \norm{\cR ( \Phi_i ) ( x_i ) - f_i ( x_i ) } _ p  ^p
        \le n ( \varepsilon n^{- \nicefrac{1}{p} } ) ^p 
        = \varepsilon ^p .
    \end{split}
\end{equation*}
In addition, for all $x = (x_1, \ldots, x_n)$, $y = (y_1, \ldots, y_n ) \in [a, b ] ^{ \sum_{i=1}^n d_i } $
with $ \forall \, i \in \num{n} \colon x_i , y_i \in [a, b ] ^{ d_i }$ we obtain that
\begin{equation*}
    \begin{split}
        \norm{\cR ( \Phi ) ( x ) - \cR ( \Phi ) ( y ) } _ p ^p 
        & = \ssuml_{i=1} ^n \norm{ \cR ( \Phi_i ) ( x_i ) - \cR ( \Phi_i ) ( y_i ) } _ p ^p
         \le L ^p \ssuml_{i=1}^n \norm{x_i - y_i } _p ^p \\
         &= L^p \norm{x-y} _p ^p 
    \end{split}
\end{equation*}
and hence $\cR ( \Phi )$ is indeed $L$-Lipschitz on $[a, b ] ^{ \sum_{i=1}^n d_i }$ with respect to $\norm{\cdot} _p$.
\end{case}
\begin{case}
Assume $p=\infty $.
Then
we get for all $x = ( x_1, \ldots, x_n ) \in [a, b ] ^{ \sum_{i=1}^n d_i } $
with $ \forall \, i \in \num{n} \colon x_i \in [a, b ] ^{ d_i }$
that
\begin{equation*}
    \begin{split}
        \norm{\cR ( x ) - ( f_1 \square f_2 \square \cdots \square f_n ) ( x ) } _\infty
        & = \max\nolimits_{i \in \num{n}} \norm{\cR ( \Phi_i ) ( x_i ) - f_i ( x_i ) } _\infty 
        \le \varepsilon .
    \end{split}
\end{equation*}
In addition, for all $x = (x_1, \ldots, x_n)$, $y = (y_1, \ldots, y_n ) \in [a, b ] ^{ \sum_{i=1}^n d_i }$
with $  \forall \, i \in \num{n} \colon x_i , y_i \in [a, b ] ^{ d_i }$ we obtain that
\begin{equation*}
    \begin{split}
        \norm{\cR ( \Phi ) ( x ) - \cR ( \Phi ) ( y ) } _ \infty
        &= \max\nolimits_{i \in \num{n} } \norm{ \cR ( \Phi_i ) ( x_i ) - \cR ( \Phi_i ) ( y_i ) } _ \infty \\
        &
        \le L \max\nolimits_{i \in \num{n}} \norm{x_i - y_i } _\infty  = L\norm{x-y} _\infty 
    \end{split}
\end{equation*}
and thus $\cR ( \Phi )$ is indeed $L$-Lipschitz with respect to $\norm{\cdot}_\infty $.
\end{case}
\end{cproof2}

\subsection{DNN approximations of parallelized Lipschitz functions}

We next state the essentially known result that every Lipschitz continuous function on a compact set in $\R^d$ can be approximated with error $\varepsilon > 0$ using at most $\cO ( \varepsilon ^{ - 2 d } )$ parameters. While this rate is not optimal (cf., e.g., Shen et al.~\cite{ShenYangZhang2020}),
it is also important for our purposes to control the Lipschitz constant of the approximating DNN realization itself, which does not directly follow from the optimal bounds in the literature.
The proof of \cref{prop:Lipschitz:approx} uses the results from \cite{JentzenRiekert2023_Error}
and \cref{cost:p:norm}.

\begin{prop}
\label{prop:Lipschitz:approx}
Let $d \in \N$.
Then there exists $K \in \N$
which satisfies for all
$a \in \R$,
$b \in (a , \infty )$,
$L \in [0 , \infty )$,
$p \in [ 1 , \infty ]$,
$\varepsilon \in (0 , 1 ]$,
and
$f \colon [a,b] ^d \to \R$ with $\forall \, x,y \in [a,b]^d \colon \abs{f ( x ) - f ( y ) } \le L \norm{x - y }_p$
that
\begin{equation}
    \cost{p} \rbr[\big]{ f, d^{ 1 - \nicefrac{1}{p } } L, \varepsilon } \le K ( \max \cu{ L ( b - a ) , 1 } ) ^{ 2 d } \varepsilon ^{ - 2 d } .
\end{equation}
\end{prop}

\begin{cproof2}{prop:Lipschitz:approx}
First \nobs that for all $ u \in \R^d$,
$p \in [ 1 , \infty ]$
we have $\norm{ u }_p \le \norm{ u }_1$.
Hence, \cite[Corollary 4.9 and Eq.~(4.38)]{JentzenRiekert2023_Error} show that there exists
$K \in \N$
which satisfies for all
$a \in \R$,
$b \in (a , \infty )$,
$L \in [0 , \infty )$,
$p \in [ 1 , \infty ]$,
$\varepsilon \in (0 , 1 ]$,
and
$f \colon [a,b] ^d \to \R$ with $\forall \, x , y \in [ a , b ] ^d \colon \abs{ f ( x ) - f ( y ) } \le L \norm{ x - y }_p \le L \norm{ x - y }_1$
that\footnote{The Lipschitz property of the approximating DNN follows since its realization is given by a maximum convolution; cf.~\cite[Lemma 3.12 and Proposition 4.4]{JentzenRiekert2023_Error}.}
\begin{equation*} 
    \cost{1} ( f , L , \varepsilon ) \le K  ( \max \cu{ L ( b - a ) , 1 } ) ^{ 2 d } \varepsilon ^{ - 2 d } .
\end{equation*}
Combining this with \cref{cost:p:norm}
demonstrates 
for all
$a \in \R$,
$b \in ( a , \infty )$,
$L \in [ 0 , \infty )$,
$p \in [ 1 , \infty ]$,
$\varepsilon \in ( 0 , 1 ]$,
and
$f \colon [ a , b ] ^d \to \R$ with $\forall \, x , y \in [ a , b ]^d \colon \abs{ f ( x ) - f ( y ) } \le L \norm{ x - y }_p$
that
\begin{equation*}
    \begin{split}
        \cost{p} ( f , d^{ 1 - \nicefrac{1}{p } } L , \varepsilon )
        \le \cost{1} ( f , L, \varepsilon ) \le K ( \max \cu{ L ( b - a ) , 1 } ) ^{ 2 d } \varepsilon ^{ - 2 d } .
    \end{split}
\end{equation*}
\end{cproof2}

Combining this result with \cref{lem:approx:parallel} we obtain the following general result about parallelized low-dimensional Lipschitz functions with input dimension $d_i \le k$.

\begin{cor}
\label{cor:lip:parallel}
Let $k \in \N$,
$L \in [ 0 , \infty )$,
$n \in \N$,
$a \in \R$,
$b \in ( a , \infty )$,
$p \in [ 1 , \infty ]$,
$d_1, \ldots, d_n \in \N$ satisfy $\max _{i \in \num{n}} d_i \le k$,
and for every $i \in \num{n}$ let
$f_i \in C ( [ a , b ] ^{ d_i } , \R )$
satisfy
$\forall \, x, y \in [ a , b ] ^{ d_i } \colon \abs{ f_i ( x ) - f_i ( y ) } \le L \norm{x - y } _p$.
Then there exists $K \in \N$,
depending only on $k$,
which satisfies for all $ \varepsilon \in ( 0 , 1 ] $
that
\begin{equation}
    \cost{p} \rbr[\big]{  f_1 \square f_2 \square \cdots \square f_n , k ^{ 1 - \nicefrac{1}{p } } L , \varepsilon }
    \le K n ^{2 k + 3 } ( \max \cu{ L ( b - a ) , 1 } ) ^{ 2 k } \varepsilon ^{ - 2 k} .
\end{equation}
\end{cor}
\begin{cproof}{cor:lip:parallel}
\Nobs that \cref{prop:Lipschitz:approx}
and \cref{lem:cost:mono} imply
that there exists $K \in \N$,
depending only on $k$,
which satisfies for all $i \in \num{n}$,
$\varepsilon \in  (0 , 1 ]$
that
\begin{equation*}
    \cost{p} \rbr[\big]{ f_i , k^{ 1 - \nicefrac{1}{p} } L , \varepsilon } \le
     \cost{p} \rbr[\big]{ f_i , ( d_i ) ^{ 1 - \nicefrac{1}{p} } L , \varepsilon } \le
      K ( \max \cu{ L ( b - a ) , 1 } ) ^{ 2 k } \varepsilon ^{ - 2 k } .
\end{equation*}
Hence, \cref{lem:approx:parallel} demonstrates for all $\varepsilon \in (0 , 1 ]$
that
\begin{equation*}
    \begin{split}
        & \cost{p} \rbr[\big]{  f_1 \square f_2 \square \cdots \square f_n , k ^{ 1 - \nicefrac{1}{p } } L , \varepsilon } \\
        & \le \tfrac{11}{4} n ^2 \rbr*{ \max \nolimits_{i \in \num{n} } d_i } ^2
        \ssuml_{i=1}^n \cost{p} \rbr[\big]{f_i, k ^{ 1 - \nicefrac{1}{p } } L , n^{- \nicefrac{1}{p} } \varepsilon } \\ 
        & \le \tfrac{11}{4} n ^2 k ^2 n K ( \max \cu{ L ( b - a ) , 1 } ) ^{ 2 k }  \varepsilon ^{ - 2 k } n ^{ \nicefrac{ 2 k }{ p } } \\
        & \le \tfrac{11 k ^{ 2 } K }{4} n^{ 2 k + 3 } ( \max \cu{ L ( b - a ) , 1 } ) ^{ 2 k }
         \varepsilon ^{ - 2 k } .
    \end{split}
\end{equation*}
\end{cproof}

Finally, we restate this result in terms of the parallelized Lipschitz function classes $\scrP_{k, p}$ as defined above. Note that the required number of parameters grows only polynomially in the input dimension $d$ for every fixed integer $k \in \N$.

\begin{cor}
	\label{cor:parallel:lip:class}
Let $k \in \N$.
Then there exists $K \in \N$
which satisfies for all
$d \in \N$,
$a \in \R$,
$b \in ( a , \infty )$,
$p \in [ 1 , \infty ]$,
$L \in [0 , \infty )$,
$\varepsilon \in (0 , 1 ]$,
$f \in \scrP_{k, p } ( [a , b ] ^d , L )$
that
\begin{equation}
    \cost{p} \rbr[\big]{ f , k ^{ 1 - \nicefrac{1}{p} } L , \varepsilon }
    \le K d^{ 2 k + 3 } ( \max \cu{ L ( b - a ) , 1 } ) ^{ 2 k } \varepsilon ^{ - 2 k} .
\end{equation}
\end{cor}

\subsection{DNN approximations of products}

In this section we state some results on the approximation cost of multidimensional products
as defined in \cref{eq:max:prod}.
The first result, \cref{lem:product:approx} below,
is a consequence of \cite[Proposition 6.8]{BeneventanoGraeber2021}
and \cref{cost:p:norm}.
Similar approximation results for multidimensional products without the curse of dimensionality can also be found in \cite{Cheridito2022_Catalog,SchwabZech2019,Yarotsky2018optimal}.

\begin{lemma}
\label{lem:product:approx}
There exists a constant
$K \in \N$
which satisfies for all
$d \in \N$,
$ a \in [1 , \infty )$,
$p \in [1 , \infty ]$,
$\varepsilon \in (0 , \infty )$
that
\begin{equation}
\label{lem:product:approx:eqclaim}
    \cost{p}  \rbr[\big]{ p_d |_{ [ - a , a ] ^d }, \sqrt{32} d^3 a^{2 d - 1 }, \varepsilon }
    \le K  d^3 \rbr*{ 1 + \ln ( a ) + \fr ( \ln ( \varepsilon ^{-1} ) ) } .
\end{equation}
\end{lemma}

\begin{cproof2}{lem:product:approx}
	\Nobs that \cite[Proposition 6.8]{BeneventanoGraeber2021} demonstrates that there exists $K \in \N$
	which satisfies for all
	$d \in \N$,
	$ a \in [1 , \infty )$,
	$\varepsilon \in (0 , \infty )$
	that
	\begin{equation*}
	   \cost{2}  \rbr[\big]{ p_d |_{ [ - a , a ] ^d }, \sqrt{32} d^{\nicefrac{5}{2} } a^{2 d - 1 }, \varepsilon }
	\le K  d^3 \rbr*{ 1 + \ln ( a ) + \fr ( \ln ( \varepsilon ^{-1} ) ) } .
	\end{equation*}
	Combining this with \cref{cost:p:norm}
	and the fact that for all
	$d \in \N$, $p \in [1 , \infty ]$
	it holds that 
	$d^{ \nicefrac{1}{2} - \nicefrac{1}{p}} \le d^{ \nicefrac{1}{2} }$
	establishes \cref{lem:product:approx:eqclaim}.
\end{cproof2}

Since the Lipschitz constant in \cref{lem:product:approx:eqclaim} grows exponentially in $d$ for $a>1$, in the following we restrict to the hypercube $[-1, 1] ^d$.
For the parallelized product functions in $\bfP ( [-1 , 1 ] ^d )$ we obtain as a consequence the following approximation result.

\begin{cor}
\label{cor:parallelized:product}
There exists a constant $K \in \N$
which satisfies for all
$d \in \N$,
$\varepsilon \in (0 , 1 ]$,
$p \in [1 , \infty ]$,
$f \in \bfP ( [ - 1 , 1 ] ^d )$
that
\begin{equation}
    \cost{p} \rbr[\big]{ f , \sqrt{32} d ^3 , \varepsilon }
    \le K d ^K \varepsilon ^{ - 1 } 
\end{equation}
\end{cor}

\begin{cproof}{cor:parallelized:product}
\Nobs that \cref{lem:product:approx} ensures that there exists a constant $K \in \N$ which satisfies for all
$d \in \N$, $p \in [1 , \infty ]$,
$\varepsilon \in (0 , \infty )$
that
\begin{equation*}
     \cost{p}  \rbr[\big]{ p_d |_{ [ - 1 , 1 ] ^d }, \sqrt{32} d^3 , \varepsilon }
    \le K  d^3 \rbr*{ 1 + \fr ( \ln ( \varepsilon ^{-1} ) ) } .
\end{equation*}
Combining this with \cref{lem:approx:parallel}
and the fact that $\forall \, u \in [1 , \infty ) \colon 1 + \ln ( u ) \le u$
demonstrates for all
$p \in [1 , \infty ]$,
$\varepsilon \in (0 , 1 ]$,
$d, n, d_1, d_2, \ldots, d_n \in \N$ with $\sum_{i=1}^n d_i = d$
and all
$f = (p_{d_1} | _ { [ - 1 , 1 ] ^{ d_1 } } ) \square \cdots \square ( p_{d_n } | _ { [ - 1 , 1 ] ^{ d_n } } ) \in  \bfP ( [ - 1 , 1 ] ^d )$
that
\begin{equation*}
    \begin{split}
        \cost{p} ( f , \sqrt{32} d ^3  , \varepsilon )
         & \le \textstyle \frac{11}{4} n ^2 \rbr*{ \max _{i \in \num{n} } d_i } ^2 \ssuml_{i=1}^n \cost{p} \rbr[\big]{ p_{d_i} | _ { [ - 1 , 1 ] ^{ d_i } } , \sqrt{32} d ^3 , n^{- \nicefrac{1}{p } } \varepsilon  } \\
         & \le \tfrac{11}{4} n ^2 d ^2 K \ssuml_{i=1}^n (d_i ) ^3 \rbr[\big]{ 1 + \ln ( \varepsilon ^{ - 1 } ) + \tfrac{1}{p} \ln n } \\
         & \le \tfrac{11}{4} n^2 d^2 K d^3 ( 1 + \ln ( \varepsilon ^{-1} ) ) ( 1 + \ln n ) \\
         & \le \tfrac{11}{4} K n^3 d^5 \varepsilon ^{-1} \le \tfrac{11}{4} K d^8 \varepsilon ^{-1} .
    \end{split}
\end{equation*}
\end{cproof}

Using the monotonicity in \cref{lem:cost:mono}, the same approximation result holds for  $f \in \bfP ( Q )$ for any $d$-dimensional hypercube $Q \subseteq [-1, 1 ] ^d $.

We next turn to the extended product functions $\fp_d$.
Combining \cite[Corollary 6.9]{BeneventanoGraeber2021} with \cref{cost:p:norm}
we obtain an approximation of $\fp_d$ with respect to arbitrary $\ell_p$-norms.

\begin{lemma}
\label{cor:prod:multi}
There exists a constant $K \in \N$
which satisfies for all $d \in \N$,
$p \in [1 , \infty ]$,
$\varepsilon \in (0 , 1 ]$
that
\begin{equation}
\label{cor:prod:multi:eqclaim}
    \cost{p}  \rbr[\big]{ \fp_d |_{ [ - 1 , 1  ] ^d }, \sqrt{32} d^{ \nicefrac{7}{2} } , \varepsilon }
    \le K  d^K \varepsilon ^{-1} .
\end{equation}
\end{lemma}

\begin{cproof}{cor:prod:multi}
	\Nobs that \cite[Corollary 6.9]{BeneventanoGraeber2021}
	implies that there exists $K \in \N$
	which satisfies for all $d \in \N$,
	$a \in [1 , \infty )$,
	$\varepsilon \in (0 , \infty )$
	that 
	\begin{equation*}
	\cost{2}  \rbr[\big]{ \fp_d |_{ [ - a , a ] ^d }, \sqrt{32} d^3 a^{2 d - 1 }, \varepsilon }
	\le K  d^5 \rbr*{ 1 + \ln ( a ) + \fr ( \ln ( \varepsilon ^{-1} ) ) } .
	\end{equation*}
	Combining this with \cref{cost:p:norm}
	yields for all
	$d \in \N$,
	$p \in [ 1 , \infty ]$,
	$\varepsilon \in (0 , \infty ) $
	that
	\begin{equation*}
		\begin{split}
		\cost{p}   \rbr[\big]{ \fp_d |_{ [ - 1 , 1  ] ^d }, \sqrt{32} d^{\nicefrac{7}{2} } , d^{ \nicefrac{1}{2} } \varepsilon } \le K  d^5 \rbr*{ 1 + \fr ( \ln ( \varepsilon ^{-1} ) ) } .
		\end{split}
	\end{equation*}
	Hence, we obtain for all
	$d \in \N$,
	$p \in [1 , \infty ]$,
	$\varepsilon \in (0 , 1 ]$
	that
	\begin{equation*}
	\begin{split}
	\cost{p} \rbr[\big]{ \fp_d |_{ [ - 1 , 1  ] ^d }, \sqrt{32} d^{\nicefrac{7}{2} } , \varepsilon }
	& \le K  d^5 \rbr*{ 1 + \fr ( \ln ( d^{ \nicefrac{1}{2} } \varepsilon ^{-1} ) ) } \\
		& \le K d^5 ( 1 + \ln ( d ) ) ( 1 + \ln ( \varepsilon ^{ - 1 } ) ) 
		\le K d^6 \varepsilon ^{-1} .
	\end{split}
	\end{equation*}
	This establishes \cref{cor:prod:multi:eqclaim}.
\end{cproof}
Again the analogous result holds for $\fp_d |_ Q $
for any hypercube $Q \subseteq [-1, 1 ] ^d $.

\subsection{DNN approximations of products with Lipschitz constant 1}

In this subsection we show that on smaller cubes of side-length at most $\frac{1}{4}$ we can even approximate the product function $p_d$ with Lipschitz constant $1$ with respect to arbitrary $\ell_p$-norms.
The proof of \cref{lem:product:approx:lip1} is based on \cite[Lemma 6.7]{BeneventanoGraeber2021}, which essentially implies the claimed statement for dimensions equal to a power of $2$. For a general input dimension we use a form of backward induction.
If we want to approximate compositions of a polynomially growing number of functions we need the Lipschitz constant to be at most 1; cf.~\cref{remark:cost:cod} above.

\begin{lemma}
\label{lem:product:approx:lip1}
There exists a constant $K \in \N$
which satisfies for all $d \in \N$,
$p \in [1 , \infty ]$,
$\varepsilon \in (0 , 1 ]$
that
\begin{equation}
    \cost{ p }  \rbr[\big]{ p_d |_{ [- \frac{1}{8} , \frac{1}{8} ] ^d  }, 1, \varepsilon }
    \le K  d^4 \rbr*{ 1 +  \ln ( \varepsilon ^{-1} ) } .
\end{equation}
\end{lemma}

\begin{cproof2}{lem:product:approx:lip1}
First, \cite[Lemma 6.7]{BeneventanoGraeber2021} implies that there exists a constant $K \in \N$
which satisfies for all $e \in \N$, $a , \varepsilon \in (0 , 1 ]$ that
\begin{equation*}
\cost{2} \rbr[\big]{ p_{ 2 ^e } |_{ [ - a , a ] ^{ 2 ^e } }, 2 ^{ 5 e / 2 } a ^{ 2 ^e - 1 } , \varepsilon } \le K 8 ^e \rbr*{ 1 + \ln ( \varepsilon ^{ - 1 } ) } .
\end{equation*}
Combining this with \cref{cost:p:norm} demonstrates for all $p \in [1 , \infty ]$,
$e \in \N$, $a , \varepsilon \in (0 , 1 ]$
that
\begin{equation*}
\cost{p} \rbr[\big]{ p_{ 2 ^e } |_{ [ - a , a ] ^{ 2 ^e } }, 2 ^{ 3 e } a ^{ 2 ^e - 1 } , \varepsilon } \le K 8 ^e \rbr*{ 1 + \ln ( \varepsilon ^{ - 1 } ) } .
\end{equation*}
Applying this with $a = \frac{1}{8}$ and $\varepsilon \with 8^{ - 2^{ e-1} } \varepsilon$ shows that for all $p \in [1 , \infty ]$,
$e \in \N$, $ \varepsilon \in (0 , 1 ]$ there exists a network $\Phi_{e}^{p, \varepsilon} \in \bfN$ such that
\begin{multline}
\label{product:approx:lip:eq:power2}
\cR (\Phi_{e}^{p, \varepsilon} ) \in C ( \R^{2^e} , \R), \qquad 
\forall \, x \in [ - \tfrac{1}{8} , \tfrac{1}{8} ] ^{2^e } \colon \abs{ \cR ( \Phi_{e}^{p, \varepsilon} ) - p_{2^e} ( x ) } \le 8^{ - 2^{ e-1} } \varepsilon, \\
\forall \, x ,y \in [ - \tfrac{1}{8} , \tfrac{1}{8} ] ^{2^e } \colon \abs{ \cR ( \Phi_{e}^{p, \varepsilon} ) ( x ) - \cR ( \Phi_{e}^{p, \varepsilon} ) ( y ) } \le 8^{ e + 1 - 2 ^e } \norm{x-y} _p, \\
\text{and} \qquad
\cP ( \Phi_{e}^{p, \varepsilon} ) \le  K 8 ^e \rbr*{ 1 + \ln ( \varepsilon ^{ - 1 } ) + 2^{e-1} \ln (8) } \le (K + \ln (8) ) 16^e  \rbr*{ 1 + \ln ( \varepsilon ^{ - 1 } ) } .
\end{multline}
Now let $d \in \N$, assume w.l.o.g.~that $d \ge 3$ (otherwise, $d$ is a power of $2$),
 and choose $e \in \N$ with $2^{ e-1} +1 \le d \le 2^e$.
Denote by $A \colon \R^d \to \R^{ 2 ^e }$ the affine linear map defined by
\begin{equation*}
A (x_1, \ldots, x_d) = \rbr*{ x_1, \ldots, x_d, \tfrac{1}{8}, \ldots, \tfrac{1}{8} } \in \R^{2 ^e } .
\end{equation*}
Note that $A$ maps $ [ - \tfrac{1}{8} , \tfrac{1}{8} ] ^{ d }$ into $  [ - \tfrac{1}{8} , \tfrac{1}{8} ] ^{2^e } $.
Given  $p \in [1 , \infty ]$, $\varepsilon \in (0 , 1 ]$ we denote by $\Psi_d ^{p, \varepsilon} \in \bfN$ the network with realization given by
\begin{equation*}
\cR ( \Psi_d ^{p, \varepsilon} ) ( x ) = 8^{ 2^e - d } \cR ( \Phi_e^{p, \varepsilon} ) ( A ( x ) ).
\end{equation*}
In other words, $\Psi_d ^{p, \varepsilon} $ is obtained from $\Phi_e^{p, \varepsilon}$ by pre-composition with the affine map $A$ and post-composition with a scalar multiplication.
Hence, \cite[Proposition 2.20]{JentzenRiekert2023_Error} implies that $\cI ( \Psi_d ^{p, \varepsilon} ) = d \le 2^e = \cI ( \Phi_e^{p, \varepsilon} )$ and all other layer dimensions are equal, whence 
\begin{equation*}
\cP ( \Psi_d ^{p, \varepsilon} ) \le  \cP ( \Phi_e^{p, \varepsilon} ) \le  (K + \ln (8) ) 16^e  \rbr*{ 1 + \ln ( \varepsilon ^{ - 1 } ) } \le 16 (K + 3 ) d ^4 \rbr*{ 1 + \ln ( \varepsilon ^{ - 1 } ) } .
\end{equation*}
 Furthermore, note that \cref{product:approx:lip:eq:power2}
 and the fact that $p_{2^e } ( A ( x ) ) = 8^{ d - 2^e } p_d ( x )$ show for all $x , y \in  [ - \frac{1}{8} ,  \frac{1}{8} ] ^d $ that
\begin{equation*}
\abs{\cR ( \Psi_d^{p, \varepsilon}  ) ( x ) - p_d ( x ) } = \abs{ 8^{ 2^e - d } \cR ( \Phi_e^{p , \varepsilon} ) ( A ( x ) ) - 8^{ 2^e - d } p_{2^e } ( A ( x ) ) } \le 8^{ 2^e - d } 8^{ - 2^{ e-1} } \varepsilon  \le \varepsilon
\end{equation*}
and
\begin{equation*}
\begin{split}
\abs{\cR ( \Psi_ d^{p, \varepsilon} ) ( x ) - \cR ( \Psi_d^{p, \varepsilon}  ( y ) ) }
& = 8^{ 2^e - d } \abs{ \cR ( \Phi_e^{p, \varepsilon} ) ( A ( x ) ) - \cR ( \Phi_e^{p, \varepsilon} ) ( A ( y ) ) } \\
& \le 8^{ 2^e - d } 8 ^{ e + 1 - 2 ^e } \norm{ A ( x ) - A ( y ) }_p = 8^{ e + 1 - d } \norm{x - y }_p \\
& \le 8^{e - 2^{e-1} }  \norm{x-y}_p \le \norm{x-y}_p .
\end{split}
\end{equation*}
\end{cproof2}

\begin{remark}
	From this proof one can see that the upper bound $\frac{1}{8}$ is not optimal for high dimensions $d$. Choosing $e$ as in the proof, $a $ can be an arbitrary number in $( 0 , 8 ^{- e / 2 ^{ e - 1 } }]$ in order for the Lipschitz constant to be at most $1$. This upper bound approaches $1$ as $d \to \infty$ (whence $e \to \infty $).
\end{remark}

For any hypercube $Q \subseteq [- \frac{1}{8} , \frac{1}{8} ] ^d $
we thus obtain as a consequence the result in \cref{cor:prod:parallel}
the parallelized product functions in $\bfP ( Q )$,
again using \cref{lem:approx:parallel}.

\begin{cor}
\label{cor:prod:parallel}
There exists a constant $K \in \N$
which satisfies for all
$d \in \N$,
every $d$-dimensional hypercube $Q \subseteq [- \frac{1}{8} , \frac{1}{8} ] ^d $,
and all
$\varepsilon \in ( 0 ,  1 ] $,
$p \in [1 , \infty ]$,
$f \in \bfP ( Q )$
that
\begin{equation}
    \operatorname{Cost}_p ( f , 1 , \varepsilon ) \le K d ^K \varepsilon ^{ - 1 } .
\end{equation}
\end{cor}
\begin{cproof}{cor:prod:parallel}
Observe that
\cref{lem:product:approx:lip1} implies that there exists a constant $K \in \N$ which satisfies for all
$d \in \N$, every hypercube $Q = [a, b ] ^d  \subseteq [- \frac{1}{8} , \frac{1}{8} ] ^d $,
and all
$\varepsilon \in (0 ,  1 ]$,
$p \in [1 , \infty ]$,
that
\begin{equation*}
     \cost{p}  \rbr[\big]{ p_d |_{Q }, 1 , \varepsilon }
    \le K  d^4 \rbr*{ 1 + \ln ( \varepsilon ^{-1} ) } .
\end{equation*}
Combining this with \cref{lem:approx:parallel}
and the fact that $\forall \, u \in [1 , \infty ) \colon 1 + \ln ( u ) \le u$
shows for all
$p \in [1 , \infty ]$,
$\varepsilon \in (0 , 1 ]$,
$d, n, d_1, d_2, \ldots, d_n \in \N$ with $\sum_{i=1}^n d_i = d$
and all
$f = (p_{d_1} | _ { [ a , b ] ^{ d_1 } } ) \square \cdots \square ( p_{d_n } | _ { [ a , b ] ^{ d_n } } ) \in  \bfP ( Q )$
that
\begin{equation*}
    \begin{split}
        \cost{p} ( f , 1 , \varepsilon )
         & \le \textstyle \frac{11}{4} n ^2 \rbr*{ \max _{i \in \num{n} } d_i } ^2 \ssuml_{i=1}^n \cost{p} \rbr[\big]{ p_{d_i} | _ { [ a , b ] ^{ d_i } } , 1 , n^{- \nicefrac{1}{p } } \varepsilon  } \\
         & \le \tfrac{11}{4} n ^2 d ^2 K \ssuml_{i=1}^n (d_i ) ^4 \rbr[\big]{ 1 + \ln ( \varepsilon ^{ - 1 } ) + \tfrac{1}{p} \ln n } \\
         & \le \tfrac{11}{4} n^2 d^2 K d^4 ( 1 + \ln ( \varepsilon ^{-1} ) ) ( 1 + \ln n ) \\
         & \le \tfrac{11}{4} K n^3 d^6 \varepsilon ^{-1} \le \tfrac{11}{4} K d^{9} \varepsilon ^{-1} .
    \end{split}
\end{equation*}
\end{cproof}

\subsection{DNN approximations of maxima}
Our final example of a particular family of functions which can be approximated by DNNs without the curse of dimensionality is given by the multidimensional maximum functions as defined in \cref{eq:max:prod}.
In fact, it is well-known in the scientific literature that these functions can be represented exactly by ReLU DNNs with only polynomially many parameters; cf., e.g.,~\cite{BeckJentzenKuckuck2022,BeneventanoGraeber2021,Cheridito2022_Catalog,JentzenRiekert2023_Error}.
Using \cref{cost:p:norm} and~\cite[Proposition 5.4]{BeneventanoGraeber2021}
we derive the following result for an approximation with Lipschitz constant $1$.

\begin{lemma}
\label{lem:max:approx}
There exists a constant $K \in \N$
which satisfies for all
$d \in \N$,
every $d$-dimensional hypercube $Q$,
and all
$\varepsilon \in [0 , \infty )$,
$p \in [1 , \infty ]$
that
\begin{equation}
\label{lem:max:approx:eqclaim}
    \cost{p} ( m_d | _ { Q } , 1 , \varepsilon ) \le K d ^2 .
\end{equation}
\end{lemma}

\begin{cproof2}{lem:max:approx}
	\Nobs that \cite[Proposition 5.4]{BeneventanoGraeber2021}
	shows that there exist
	$K \in \N$ and $\Phi_d \in \bfN$, $d \in \N$,
	which satisfy for all
	$d \in \N$,
	$x, y \in \R^d$
	that
	$\cR ( \Phi_d ) \in C ( \R^d , \R )$,
	$\cP ( \Phi_d ) \le K d ^2 $,
	$\cR ( \Phi_d ) ( x ) = m_d ( x ) $,
	and
	$\abs{ \cR ( \Phi_d ) ( x ) - \cR ( \Phi_d ) ( y ) } \le \norm{ x - y } _ \infty $.
	This ensures for all
	$d \in \N$,
	$\varepsilon \in [0 , \infty )$
	and every $Q = [a , b ] ^d \subseteq \R^d$
	that
	$\cost{\infty} ( m_d | _ { Q } , 1 , \varepsilon ) \le K d ^2$.
	Combining this with \cref{cost:p:norm} establishes \cref{lem:max:approx:eqclaim}.
\end{cproof2}

Similarly as for the product functions, we get for the parallelized maximum functions in $\bfM ( Q )$ the following approximation result.

\begin{cor}
\label{cor:max:parallel}
There exists a constant $K \in \N$
which satisfies for all
$d \in \N$,
every $d$-dimensional hypercube $Q \subseteq \R ^d $,
and all
$ \varepsilon \in [0 , \infty )$,
$p \in [1 , \infty ]$,
$f \in \bfM ( Q )$
that
\begin{equation}
    \cost{p} \rbr*{ f , 1 , \varepsilon }
    \le K d ^K .
\end{equation}
\end{cor}

\begin{cproof2}{cor:max:parallel}
\Nobs that \cref{lem:max:approx} assures that there exists a constant $K \in \N$
which satisfies for all
$d \in \N$,
$Q = [a, b ] ^d \subseteq \R^d $,
$\varepsilon \in [0 , \infty )$,
$p \in [1 , \infty ]$
that
\begin{equation*} 
    \cost{p} ( m_d | _ { Q } , 1 , \varepsilon ) \le K d ^2 .
\end{equation*}
Analogously to the proof of \cref{cor:prod:parallel}, we combine this with \cref{lem:approx:parallel}
to obtain for all
$a \in \R$,
$b \in (a , \infty )$,
$p \in [1 , \infty ]$,
$\varepsilon \in (0 , \infty )$,
$d, n, d_1, d_2, \ldots, d_n \in \N$ with $\sum_{i=1}^n d_i = d$
and all
$f = (m_{d_1} | _ { [ a , b ] ^{ d_1 } } ) \square \cdots \allowbreak \square ( m_{d_n } | _ { [a , b ] ^{ d_n } } ) \in  \bfM ( Q )$
that
\begin{equation*}
    \begin{split}
        \cost{p} ( f , 1 , \varepsilon )
        & \le \tfrac{11}{4} n ^2 \rbr*{ \max\nolimits_{i \in \num{n} } d_i } ^2 \ssuml_{i=1}^n \cost{p} \rbr[\big]{ m_{d_i } | _ { [a , b ] ^{ d_i } } , 1 , n ^{ - \nicefrac{1}{p } } \varepsilon } \\
        & \le \tfrac{11}{4} n ^2 d ^2 \ssuml_{i=1}^n K (d_i ) ^2 \le \tfrac{11}{4} K d^6 .
    \end{split}
\end{equation*}
\end{cproof2}

For the extended maximum functions $\fm_d$
we obtain the following result by using \cite[Corollary 5.5]{BeneventanoGraeber2021} and \cref{cost:p:norm}.
For $p < \infty$ the Lipschitz constant depends on the dimension, due to the use of the $\ell_p$-norm on the output space $\R^d$.

\begin{lemma}
\label{lem:max:multi}
There exists a constant $K \in \N$
which satisfies for all
$d \in \N$,
every $d$-dimensional hypercube $Q \subseteq \R ^d $,
and all
$\varepsilon \in [0 , \infty )$,
$p \in [1 , \infty ]$
that
\begin{equation}
\label{lem:max:multi:eqclaim}
     \cost{p} \rbr[\big]{ \fm_d |_{ Q }, d ^{ \nicefrac{1}{p} } , \varepsilon }
    \le K  d^4 .
\end{equation}
\end{lemma}

\begin{cproof}{lem:max:multi}
	\Nobs that \cite[Corollary 5.5]{BeneventanoGraeber2021}
	assures that there exist $K \in \N$
	and $\Phi_d \in \bfN$, $d \in \N$,
	which satisfy for all
	$d \in \N$,
	$x \in \R^d$
	that
	$\cR ( \Phi_d ) \in C ( \R^d , \R^d )$,
	$\cP ( \Phi_d ) \le K d ^4 $,
	and
	$ \cR ( \Phi_d ) ( x ) = \fm_d ( x ) $.
	Furthermore, \nobs that $\fm_d$ is $1$-Lipschitz with respect to the $\ell_\infty$-norm.
	Hence, we obtain for all
	$d \in \N$,
	$Q = [a, b ] ^d \subseteq \R^d $,
	$\varepsilon \in [0 , \infty )$
	that
	$  \cost{ \infty } \rbr[\big]{ \fm_d |_{ Q }, 1 , \varepsilon }
	\le K  d^4 $.
	Combining this with \cref{cost:p:norm}
	establishes \cref{lem:max:multi:eqclaim}.
\end{cproof}

\subsection{Proof of the main results}
\label{subsec:proof:main}
Now we have all the ingredients to prove the main theorems from the introduction.
We first employ \cref{prop:main:abstract:comp1} to establish the result in \cref{theo:intro:1} with a polynomial number of functions in the composition.

\begin{cproof}{theo:intro:1}
\Nobs that
\cref{cor:parallel:lip:class}
shows that there exists $K_1 \in \N$ which satisfies for all
$d \in \N$,
$i \in \num{\bfk ( d ) }$ with $g_i^d \in \scrP_{c, 1 }( Q_i^d , 1 )$
that
\begin{equation*}
\forall \, \varepsilon \in (0 , 1 ] \colon \cost{1}  ( g_i^d , 1 , \varepsilon ) 
\le K_1 (\fd_i^d ) ^{ 2 c + 3 } ( 2 c d ^c ) ^{2 c } \varepsilon^{-2 c }
\le K_1 c^{ 2 c + 3 }  (2c ) ^{ 2 c } d^{ c ( 2c + 3 ) + 2 c ^2 } \varepsilon^{-2 c }.
\end{equation*}
Here we used that \cref{cor:parallel:lip:class} can be applied with the side-length $b - a \le 2 c d ^c $, since $Q_i^d \subseteq [- c d^c , c d^c ] ^{ \fd_i^d } $ by assumption.
Furthermore,
\cref{cor:prod:parallel}
ensures that there exists $K_2 \in \N$ which satisfies for all
$d \in \N$,
$i \in \num{\bfk ( d ) }$ with $g_i^d \in \bfP ( Q_i^d )$
that
\begin{equation*}
\forall \, \varepsilon \in (0 , 1 ] \colon \cost{1}  ( g_i^d , 1 , \varepsilon ) 
\le K_2 (\fd_i^d ) ^{ K_2 } \varepsilon^{- 1 }
\le K_2 c^{ K_2} d^{ c K_2 } \varepsilon^{-2 c }.
\end{equation*}
Finally,
\cref{cor:max:parallel}
implies that there exists $K_3 \in \N$ 
which satisfies for all
$d \in \N$,
$i \in \num{\bfk ( d ) }$
with $g_i^d \in \bfM ( Q_i^d )$
that
\begin{equation*}
\forall \, \varepsilon \in (0 , 1 ] \colon \cost{1}  ( g_i^d , 1 , \varepsilon ) 
\le K_3 (\fd_i^d ) ^{ K_3 }
\le K_3 c^{ K_3} d^{ c K_3 } \varepsilon^{-2 c }.
\end{equation*}
Hence, we obtain for all
$d \in \N$,
$i \in \num{\bfk ( d ) }$,
$\varepsilon \in (0 , 1 ]$
that
\begin{equation*}
	    \cost{1} ( g_i^d , 1 , \varepsilon ) 
	    \le \max \cu*{ 2^{ 2 c } K_1 c ^{ 4 c + 3 } , K_2 c^{ K_2} , K_3 c^{ K_3} } 
	    d ^{ c \max \cu{ 4c + 3 , K_2, K_3 } } \varepsilon ^{ - 2 c } .
\end{equation*}
Combining this with \cref{prop:main:abstract:comp1} establishes the claim.
\end{cproof}

Next we apply \cref{prop:main:abstract:comp2} to establish \cref{theo:intro:2}, where the Lipschitz constants of the composed functions grow polynomially in the dimension.

\begin{cproof}{theo:intro:2}
Analogously to the proof of \cref{theo:intro:1}
we can combine
\cref{theo:intro:2:eq:allowed:function},
\cref{cor:parallel:lip:class},
\cref{cor:parallelized:product},
\cref{cor:prod:multi},
\cref{cor:max:parallel},
and
\cref{lem:max:multi}
to obtain that there exists $K \in \N$ 
which satisfies for all
$d \in \N$,
$i \in \num{\bfk ( d ) }$,
$\varepsilon \in (0 , 1 ]$
that
\begin{equation*}
\cost{p} ( g_i^d , K d ^K , \varepsilon ) \le K (\fd_i^d ) ^K \max \cu{\varepsilon ^{ - 2c } , \varepsilon ^{-1} , 1 } \le K d^K \varepsilon ^{ - 2 c } .
\end{equation*}
Applying the second abstract approximation result in
\cref{prop:main:abstract:comp2} hence establishes the claim.
\end{cproof}

\subsection*{Acknowledgments}
This work has been funded by the Deutsche Forschungsgemeinschaft (DFG, German Research Foundation) under Germany’s Excellence Strategy EXC 2044-390685587, Mathematics Münster: Dynamics-Geometry-Structure. Helpful suggestions by Patrick Cheridito, Arnulf Jentzen, Benno Kuckuck, and Florian Rossmannek, in particular regarding \cref{lem:product:approx:lip1}, are gratefully acknowledged.


\end{document}